\documentclass[12pt]{amsart}
\usepackage{amssymb,amsmath,graphicx,amsfonts,euscript}
\usepackage{amsmath,amssymb,times,color}
\usepackage{graphicx}
\usepackage{dsfont}
\usepackage[T1]{fontenc}
\definecolor{webred}{rgb}{0.75,0,0}
\definecolor{webgreen}{rgb}{0,0.75,0}
\usepackage[citecolor=webgreen,colorlinks=true,linkcolor=webred]{hyperref}

\numberwithin{equation}{section}

\makeindex \makeatletter
\def\captionof#1#2{{\def\@captype{#1}#2}}
\makeatother
\newtheorem{proposition}{Proposition}[section]
\newtheorem{remark}{\bf Remark}
\newtheorem{theorem}{\bf Theorem}[section]
\def\bgneqnn{\begin{equation}}
\def\endeqnn{\end{equation}}
\def\bgneqy{\begin{eqnarray}}
\def\endeqy{\end{eqnarray}}
\def\bgneqy*{\begin{eqnarray*}}
\def\endeqy*{\end{eqnarray*}}

\def\text{\mbox}

\def\bgneqy*{\begin{eqnarray*}}
\def\endeqy*{\end{eqnarray*}}

\def\qed{\hfill$\bullet$\par\bigskip}
\newenvironment{preuve}[1]{{\bf Proof #1.~}}{\qed}

\newcounter{tablegroup}
\newcounter{subtable}[tablegroup]

\newcommand{\handletables}

\title[NFVM for rotating channel flows]{New Finite Volume Method for rotating channel flows involving boundary layers}

\subjclass[2000]{\ 76D10, 76D17, 65L11, 76L05, 68U120.}
\keywords{Boundary layer, colocated scheme, correctors, Navier-Stokes equations, finite volumes, singular
problems.}

\author[Ben Chaabane$^{2}$, Hamouda$^{1,3}$ and  Tekitek$^3$]{\sc {Soumaya Ben Chaabane}$^{2}$, {Makram Hamouda}$^{1,3}$ and {Mahdi Tekitek}$^3$}

\address{$^1$ The Institute for Scientific Computing and
Applied Mathematics,, Indiana University, 831 E. 3rd St., Rawles Hall,, Bloomington, IN 47405, USA.}

\address{$^2$ University of Carthage, Faculty of Sciences of Bizerte, Department of Mathematics, Zarzouna,  Tunisia.}

\address{$^3$ University of Tunis El Manar, Faculty of Sciences of Tunis, Department of Mathematics, Tunis,  Tunisia.}

\begin{document}

\maketitle

\begin{abstract}
We investigate in this article the boundary layers appearing for a fluid under moderate rotation when the viscosity is small. The fluid is modeled by rotating type Stokes equations known also as the Barotropric mode equations in the primitive equations theory. First we derive the correctors that describe the sharp variations at large Reynolds number (i.e., small viscosity). Second, thanks to a new finite volume method (NFVM) we give numerical solutions of the rotating Stokes system at small viscosity. The NFVM can be applied to a large class of singular perturbation problems.
\end{abstract}

\section{Introduction}
We are interested in this article in the study of boundary layers of a time dependant rotating fluid when the viscosity is small and the boundary is characteristic; this occurs for example when the boundary is solid
and at rest.  The boundary conditions are considered of homogeneous Dirichlet type.
More precisely, we
consider the flow in $3D$ verifying the following system:
\begin{equation}\label{fluide1}
\left\lbrace
\begin{array}{lcl}
\dfrac{\partial \boldsymbol{u}^{\varepsilon}}{\partial t}-\varepsilon\Delta \boldsymbol{u}^{\varepsilon}+\boldsymbol{\omega}\times \boldsymbol{u}^{\varepsilon}
+\nabla p^{\varepsilon}=\boldsymbol{f}, \hspace{0,5cm}\mbox{in}\hspace{0,2cm}\Omega_{\infty}\times(0,T), \vspace{0,2cm}\\
\textnormal{div}~\boldsymbol{u}^{\varepsilon}=0,\hspace{0.3cm}\text{in}~\Omega_{\infty}\times(0,T),\vspace{0.2cm}\\
\boldsymbol{u}^{\varepsilon}=0,\hspace{0.3cm}\text{on}~\partial\Omega_{\infty}\\
\boldsymbol{u}^{\varepsilon}~\text{is}~2\pi\text{-periodic in the}~ x~
\text{and}~ y~\text{directions},\\
\boldsymbol{u}^{\varepsilon}_{|t=0}=\boldsymbol{u}_{0};\hspace{0,5cm}
\end{array}\right.
\end{equation}
see
\cite{3} and \cite{8} for more details about the theory of rotating fluids.
Here $\boldsymbol{\omega}=\alpha \boldsymbol{e_{3}}$ where $\boldsymbol{e_{3}}$ is the unit vector in the canonical basis of $\mathbb{R}^3, \Omega_{\infty}=\mathbb{R}^2\times(0,h)$ is the relevant
domain, $\Gamma=\partial\Omega_{\infty}=\mathbb{R}^2\times\{0,h\}$
its boundary. The functions $\boldsymbol{{u}_{0}}$ and $\boldsymbol{f}$ are given and supposed
to be as regular as necessary. Without loss of generality, the constant $h$ will be taken from now equal to $1$. \\The solution
$(\boldsymbol{u}^{\varepsilon},p^{\varepsilon})$ of the system $(\ref{fluide1})$ is such that
$\boldsymbol {u}^{\varepsilon}(t,x,y,z)=(u_{1}^{\varepsilon},u_{2}^{\varepsilon},u_{3}^{\varepsilon
})\in \mathds{R}^{3}$ and $p^{\varepsilon}\in \mathds{R}$, the
coefficient $\varepsilon$ is a positive constant
representing the inverse of the Reynolds number. Throughout this
paper the coefficient $\varepsilon>0$ is intended to be small
$\varepsilon \ll 1$. Because of the periodicity conditions
$(\ref{fluide1})_{4}$ we will consider a portion of the channel
$\Omega_{\infty}$ that we denote by
$\Omega=(0,2\pi)\times(0,2\pi)\times(0,1)$ and its boundary
$\Omega=(0,2\pi)\times(0,2\pi)\times\{0,1\}$ on which all our
calculations will be done. The formal limit solution $\boldsymbol{u}^{0}$ of the system $(\ref{fluide1})$ is simply obtained by setting $\varepsilon=0$ in $(\ref{fluide1})$. Hence, we have:
\begin{equation}\label{eq lim 00}
\left\lbrace
\begin{array}{lcl}
\dfrac{\partial \boldsymbol{u}^{0}}{\partial t}+\boldsymbol{\omega}\times \boldsymbol{u}^{0}+\nabla p^{0}=\boldsymbol{f}, \hspace{0,5cm}\mbox{in}\hspace{0,2cm}\Omega\times (0,T), \vspace{0,2cm}\\
\textnormal{div}~\boldsymbol{u}^{0}=0,\hspace{0.3cm}\text{in}~\Omega\times(0,T),\vspace{0.3cm}\\
u^{0}_{3}=0,\hspace{0.3cm}\text{on}~\partial\Omega,\vspace{0.3cm}\\
\boldsymbol{u}^{0}~\text{is}~2\pi\text{-periodic in the}~ x~ \text{and}~ y
~\text{directions},\\
\boldsymbol{u}^{0}_{|t=0}=\boldsymbol{u}_{0}.
\end{array}\right.
\end{equation}
The absence in the limit system of the Laplacian term
$(-\varepsilon\Delta \boldsymbol{u}^{\varepsilon})$ which is a regularizing
term, generates a loss of regularity of the limit solution
$\boldsymbol{u}^{0}$. Thus some discrepancies between the viscous and inviscid
solutions appear near the boundary of the domain, that is
here $z=0, 1$. This thin region is called
\emph{boundary layer} and the convergence of $\boldsymbol{u}^{\varepsilon}$ to $\boldsymbol{u}^{0}$ is
not expected there especial when we look for the convergence in $H_{X}^{s}(\Omega)$ for $s>1$ and $X=$. Hence, we introduce a correcting term called \emph{corrector}
for which the equation must be of course simpler
than the one in the original problem namely $(\ref{fluide1})_{1}$; see \cite{4}, \cite{10n}, \cite{11},
\cite{12} and \cite{13} for more details on this notion.

\section{The corrector equations}
To study the asymptotic behavior of $\boldsymbol{u}^{\varepsilon}$, when
$\varepsilon\rightarrow0$, we propose the following asymptotic
expansion of $\boldsymbol{u}^{\varepsilon}$:\\
$$\boldsymbol{u}^{\varepsilon}\simeq \boldsymbol{u}^{0}+\boldsymbol{\varphi}^{\varepsilon},$$
where $\boldsymbol{\varphi}^{\varepsilon}$ is the corrector function that will be
introduced to correct the difference $\boldsymbol{u}^{\varepsilon}-\boldsymbol{u}^{0}$ at
$z=0,1$. The equations verified by $\boldsymbol{\varphi}^{\varepsilon}$ are as follows:
\begin{equation}\label{fluide 3}
\left\lbrace
\begin{array}{lcl}
\dfrac{\partial \boldsymbol{\varphi}^{\varepsilon}}{\partial t}-\varepsilon\dfrac{\partial^{2}\boldsymbol{\varphi}^{\varepsilon}}{\partial z^{2}}
+\boldsymbol{\omega}\times\boldsymbol{\varphi}^{\varepsilon}=0,\hspace{0,5cm}\mbox{in $\Omega\times (0,T)$}, \vspace{0,2cm}\\
\textnormal{div}\boldsymbol{\varphi}^{\varepsilon}=0,\hspace{0.3cm}\text{in}~\Omega\times(0,T),\vspace{0.3cm}\\
\boldsymbol{\varphi}^{\varepsilon}_{|z=0,1}=-\boldsymbol{u}^{0}_{|z=0,1},\vspace{0,2cm}\\
\boldsymbol{\varphi}^{\varepsilon}_{|t=0}=0.
\end{array}\right.
\end{equation}

We now introduce an approximate function $\check{{\boldsymbol{\varphi}}}^{\varepsilon}$ of $\boldsymbol{\varphi}^{\varepsilon}$ defined as the sum of
$\overline{\boldsymbol{\varphi}}^{0,\varepsilon}$ and
$\widetilde{\boldsymbol{\varphi}}^{1,\varepsilon}$ the correctors that we propose to solve the boundary layers at the boundaries $z=0$ and $z=1$, respectively:

$$\check{\boldsymbol{\varphi}}^{\varepsilon}(t,x,y,z)=\overline{\boldsymbol{\varphi}}^{0,\varepsilon}(t,x,y,\frac{z}{\sqrt{\varepsilon}})+
\widetilde{\boldsymbol{\varphi}}^{1,\varepsilon}(t,x,y,\frac{1-z}{\sqrt{\varepsilon}}),$$
where $\overline{z}=\frac{z}{\sqrt{\varepsilon}}$ and
$\widetilde{z}=\frac{1-z}{\sqrt{\varepsilon}}$. Taking into
consideration the linearity of the equations $(\ref{fluide 3})_{1,2}$
and the boundary conditions $(\ref{fluide 3})_{3}$, the system verified by
$\overline{\boldsymbol{\varphi}}^{0,\varepsilon}$ is given by:
\begin{equation}\label{fluide 4}
\left\lbrace
\begin{array}{lcl}
\dfrac{\partial\overline{\boldsymbol{\varphi}}^{0,\varepsilon}}{\partial t}-\dfrac{\partial^{2}\overline{\boldsymbol{\varphi}}^{0,\varepsilon}}{\partial \overline{z}^{2}}
+\boldsymbol{\omega}\times\overline{\boldsymbol{\varphi}}^{0,\varepsilon}=0, \hspace{0,5cm}\mbox{in}\hspace{0,2cm}\widetilde{\Omega}\times(0,T), \vspace{0,2cm}\\
\overline{\boldsymbol{\varphi}}^{0,\varepsilon}(\overline{z}=0)=-\boldsymbol{u}^{0}(\overline{z}=0),\vspace{0,2cm}\\
\overline{\boldsymbol{\varphi}}^{0,\varepsilon}\rightarrow0 \hspace{0.3cm} \text{as}~\overline{z}\rightarrow\infty,\vspace{0,2cm}\\
\overline{\boldsymbol{\varphi}}^{0,\varepsilon}_{|t=0}=0.
\end{array}\right.
\end{equation}
Similarly $\widetilde{\varphi}^{1,\varepsilon}$ satisfies the
following system:
\begin{equation}\label{fluide 5}
\left\lbrace
\begin{array}{lcl}
\dfrac{\partial\widetilde{\boldsymbol{\varphi}}^{1,\varepsilon}}{\partial t}
-\dfrac{\partial^{2}\widetilde{\boldsymbol{\varphi}}^{1,\varepsilon}}{\partial \widetilde{z}^{2}}+\boldsymbol{\omega}\times\widetilde{\boldsymbol{\varphi}}^{1,\varepsilon}=0,
\hspace{0,5cm}\mbox{in}\hspace{0,2cm}\widetilde{\Omega}\times(0,T), \vspace{0,2cm}\\
\widetilde{\boldsymbol{\varphi}}^{1,\varepsilon}(\widetilde{z}=0)=-\boldsymbol{u}^{0}(\widetilde{z}=0),\vspace{0,2cm}\\
\widetilde{\boldsymbol{\varphi}}^{1,\varepsilon}\rightarrow0 \hspace{0.3cm} \text{as}~\widetilde{z}\rightarrow\infty,\vspace{0,2cm}\\
\widetilde{\boldsymbol{\varphi}}^{1,\varepsilon}_{|t=0}=0.
\end{array}\right.
\end{equation}
Here we denoted by $\widetilde{\Omega}$ the stretched domain, i.e., $\widetilde{\Omega}=(0,2\pi)\times (0,2\pi)\times (0,+\infty)$ of the system $(\ref{fluide
4})$. In the following we will derive the expressions of the solutions of the systems $(\ref{fluide 4})$ and $(\ref{fluide 5})$. For that purpose we need the following proposition.
\begin{proposition}\label{proposition}
Let $u=u(t,x,y,z)$ be the solution of the  following problem:

\begin{equation}\label{tounant}
\left\lbrace
\begin{array}{lcl}
\dfrac{\partial\boldsymbol{ u}}{\partial t}-\dfrac{\partial^{2}\boldsymbol{u}}{\partial z^{2}}+\boldsymbol{\omega}\times \boldsymbol{u}=0, \hspace{0,2cm}~\text{in}~\widetilde{\Omega}\times (0,T),\vspace{0,2cm}\\
\boldsymbol{u}=\boldsymbol{g},\hspace{0,5cm}\mbox{at}~z=0,\vspace{0,2cm}\\
\boldsymbol{u}\rightarrow 0,\hspace{0,5cm}\mbox{as}~z\rightarrow+\infty, \vspace{0,2cm}\\
\boldsymbol{u}=0,~~~\text{at}~t=0.
\end{array}\right.
\end{equation}

where $\boldsymbol{g}=(g_{1},g_{2},0)$ is a continuous function in $\widetilde{\Omega}\times(0,T)$ and $\boldsymbol{w}=\alpha \boldsymbol{e_{3}}$.
Then, the explicit expression of $\boldsymbol{u}$ is given by:
\begin{eqnarray*}
\boldsymbol{u}(t,x,y,z)&=&-\int_{0}^{t}\frac{\partial K}{\partial z}(t-\tau,z)[(\boldsymbol{g}-i(\boldsymbol{e_{3}}\times \boldsymbol{g}))(\tau,x,y,0) e^{i\alpha(\tau-t)}\\
& &+(\boldsymbol{g}+i(\boldsymbol{e_{3}}\times \boldsymbol{g}))(\tau,x,y,0)
e^{i\alpha(t-\tau)}]d\tau,
\end{eqnarray*}
where $i$ is the complex number s.t. $i^{2}=-1$, and $K$ is the fundamental solution of the heat equation:
$$\displaystyle K(t,z)=\frac{1}{\sqrt{4\pi t}}e^{\frac{-z^2}{4t}}.$$
\end{proposition}
\begin{preuve}
\normalfont Let $\boldsymbol{u}=(u_{1},u_{2},u_{3})$ the solution of $(\ref{tounant})$.
We have $g_{3}=0$, hence $\boldsymbol{u}=(u_{1},u_{2},0)$ i.e. $u_{3}=0$. Taking the cross product of $(\ref{tounant})$
with $\boldsymbol{e_{3}}$, we find:
$$\partial_{t}(\boldsymbol{e_{3}}\times \boldsymbol{u})-\partial_{z}^{2} (\boldsymbol{e_{3}}\times
u)-\alpha\boldsymbol{ u}=0.$$
We then set $\boldsymbol{C}^{\pm}=\boldsymbol{u}\mp i(\boldsymbol{e_{3}}\times
\boldsymbol{u})$, we obtain:
$$\partial_{t} \boldsymbol{C}^{\pm}-\partial_{z}^{2}C^{\pm}\pm i\alpha \boldsymbol{C}^{\pm}=0.$$
Denoting by $\boldsymbol{H}^{\pm}=\boldsymbol{C}^{\pm}e^{\pm i\alpha t}$, one arrives to the following system:
\begin{equation}\label{equation}
\left\lbrace
\begin{array}{lcl}
\dfrac{\partial \boldsymbol{H}^{\pm}}{\partial t}-\dfrac{\partial^{2}\boldsymbol{H}^{\pm}}{\partial z^{2}}=0, \hspace{0,2cm}~\text{in}~\widetilde{\Omega}\times (0,T),\vspace{0,2cm}\\
\boldsymbol{H}^{\pm}(z=0)=(\boldsymbol{g}(z=0)\mp i(\boldsymbol{e_{3}}\times \boldsymbol{g}(z=0))e^{\pm i\alpha t},\vspace{0,2cm}\\
\boldsymbol{H}^{\pm}\rightarrow 0 ,~~\text{as}~~z\rightarrow +\infty, \vspace{0,2cm}\\
\boldsymbol{H}^{\pm}|_{t=0}=0.
\end{array}\right.
\end{equation}
Hence $\boldsymbol{H}^{\pm}$ satisfies a heat equation with non-homogeneous
boundary conditions, then it has the following expression (\cite{2}):
$$\boldsymbol{H}^{\pm}=-2\int_{0}^{t}\frac{\partial K}{\partial z}(t-\tau,z)[(\boldsymbol{g}\mp i(\boldsymbol{e_{3}}\times \boldsymbol{g}))(\tau,x,y,0)]e^{\pm i\alpha\tau}d\tau.$$
Then, we infer that:
$$\boldsymbol{C}^{\pm}=-2\int_{0}^{t}\frac{\partial K}{\partial z}(t-\tau,z)[(\boldsymbol{g}\mp i(\boldsymbol{e_{3}}\times \boldsymbol{g}))(\tau,x,y,0)]e^{\pm i\alpha(\tau-t)}d\tau.$$
Coming back to $\boldsymbol{u}$ we have:
 $$\boldsymbol{u}=\frac{1}{2}(\boldsymbol{C}^{+}+\boldsymbol{C}^{-}),$$
hence we deduce the explicit expression of the solution of $(\ref{tounant})$:
\begin{eqnarray*}
\boldsymbol{u}&=&-\int_{0}^{t}\frac{\partial K}{\partial z}(t-\tau,z)\times\{[(\boldsymbol{g}- i(\boldsymbol{e_{3}}\times \boldsymbol{g}))(\tau,x,y,0)]e^{i\alpha(\tau-t)}+\\
& &+[(\boldsymbol{g}+i(\boldsymbol{e_{3}}\times
\boldsymbol{g}))(\tau,x,y,0)]e^{i\alpha(t-\tau)}\}d\tau.
\end{eqnarray*}
\end{preuve}
Now, according to Proposition \ref{proposition}, the solution of $(\ref{fluide 4})~\overline{\boldsymbol{\varphi}}^{0,\varepsilon}
=(\overline{{\varphi}}_{1}^{0,\varepsilon},\overline{{\varphi}}_{2}^{0,\varepsilon},\overline{{\varphi}}_{3}^{0,\varepsilon})
$
has the following expression:
\begin{eqnarray*}
\overline{\varphi}_{j} ^{0,\varepsilon} &=& -\int_{0}^{t}\frac{1}{\sqrt{4\pi(t-\tau)}} \frac{z}{2\sqrt{\varepsilon}(t-\tau)} e^{\frac{-z^2}{4\varepsilon (t-\tau)}}
\times\{2 u_{j}^{0}(\tau,x,y,0) \cos(\alpha(\tau-t))\\
 & &+2 (\boldsymbol{e_{3}}\times \boldsymbol{u}^{0})_{j}(\tau,x,y,0) \sin(\alpha(\tau-t))\}d\tau, \hspace{0.2cm} j=1,2,
\end{eqnarray*}
for the two tangential components of $\overline{\boldsymbol{\varphi}}^{0,\varepsilon}$, and the normal component of $\overline{\boldsymbol{\varphi}}^{0,\varepsilon}$ is simply deduced using the incompressibility condition:

\begin{eqnarray}\label{div}
\overline{\varphi}_{3} ^{0,\varepsilon} &=&
-\int_{0}^{t}\frac{\sqrt{\varepsilon}}{\sqrt{4\pi(t-\tau)}}e^{\frac{-z^2}{4\varepsilon
(t-\tau)}}
\times\{-2\partial_{z}u^{0}_{3}(\tau,x,y,0)\cos(\alpha(\tau-t))\\& &-2(\partial_{x}u_{2}^{0}-\partial_{y}u_{1}^{0})(\tau,x,y,0)\sin(\alpha(\tau-t))\}d\tau
\nonumber \\& &+\int_{0}^{t}\frac{\sqrt{\varepsilon}}{\sqrt{4\pi(t-\tau)}}e^{\frac{-1}{4\varepsilon (t-\tau)}} \times\{
-2\partial_{z}u^{0}_{3}(\tau,x,y,0)\cos(\alpha(\tau-t))\nonumber\\& &-2(\partial_{x}u_{2}^{0}-\partial_{y}u_{1}^{0})(\tau,x,y,0)\sin(\alpha(\tau-t))\}d\tau.\nonumber
\end{eqnarray}

Then we write the system satisfied by
$\overline{\boldsymbol{\varphi}}^{0,\varepsilon}$ which reads as follows:
\begin{equation}\label{fluide44}
\left\lbrace
\begin{array}{lcl}
\dfrac{\partial\overline{\boldsymbol{\varphi}}^{0,\varepsilon}}{\partial t}-\varepsilon\dfrac{\partial^{2}\overline{\boldsymbol{\varphi}}^{0,\varepsilon}}{\partial {z}^{2}}+\boldsymbol{\omega}\times\overline{\boldsymbol{\varphi}}^{0,\varepsilon}=(0,0,\dfrac{\partial\overline{\varphi}_{3} ^{0,\varepsilon}}{\partial t}-\varepsilon\dfrac{\partial^{2}\overline{\varphi}_{3}^{0,\varepsilon}}{\partial z^{2}}),\hspace{0,2cm}~\text{in}~\widetilde{\Omega}\times (0,T),\vspace{0,2cm}\\
\textnormal{div}~\overline{\boldsymbol{\varphi}}^{0,\varepsilon}=0,\hspace{0,2cm}~\text{in}~\widetilde{\Omega}\times (0,T),\vspace{0,2cm}\\
\overline{\boldsymbol{\varphi}}^{0,\varepsilon}(z=0)=(-u_{1}^{0}(z=0),-u_{2}^{0}(z=0),\overline{\varphi}_{3}^{0,\varepsilon}(z=0)),\vspace{0,2cm}\\
\overline{\boldsymbol{\varphi}}^{0,\varepsilon}(z=1)=(\overline{\varphi}_{1}^{0,\varepsilon}(z=1),\overline{\varphi}_{2}^{0,\varepsilon }(z=1),0),\vspace{0,2cm}\\
\overline{\boldsymbol{\varphi}}^{0,\varepsilon}(t=0)=0.
\end{array}\right.
\end{equation}
Now, we need to estimate the right-hand side (denoted hereafter RHS) of
$(\ref{fluide44})_{1}$. First we set the change of variables $s=\frac{1}{\sqrt{t-\tau}}$, and
we infer that:
\begin{eqnarray}\label{verif}
\overline{\varphi}_{3} ^{0,\varepsilon}&=&-\int_{\frac{1}{\sqrt{t}}}^{\infty}\frac{\sqrt{\varepsilon}}{\sqrt{\pi} s^2} e^{\frac{-z^{2}s^{2}}{4\varepsilon}}
\times\{-2\partial_{z}u^{0}_{3}(t-\frac{1}{s^2},x,y,0)\cos\left(\frac{\alpha}{s^2}\right)\\& &+2(\partial_{x}u_{2}^{0}-\partial_{y}u_{1}^{0})
(t-\frac{1}{s^2},x,y,0)\sin\left(\frac{\alpha}{s^2}\right)\}ds \nonumber
\\& &+\int_{\frac{1}{\sqrt{t}}}^{\infty}\frac{\sqrt{\varepsilon}}{\sqrt{\pi} s^2}e^{\frac{-s^{2}}{4\varepsilon}}
\times\{-2\partial_{z}u^{0}_{3}(t-\frac{1}{s^2},x,y,0)\cos\left(\frac{\alpha}{s^2}\right) \nonumber\\& &+2(\partial_{x}u_{2}^{0}-\partial_{y}u_{1}^{0})
(t-\frac{1}{s^2},x,y,0)\sin\left(\frac{\alpha}{s^2}\right)\}ds. \nonumber
\end{eqnarray}

By differentiating $(\ref{verif})$ with
respect to the time variable $t$,
 we obtain:
\begin{eqnarray}\label{new1}
\frac{\partial\overline{\varphi}_{3} ^{0,\varepsilon}}{\partial
t}&=&-\frac{\sqrt{\varepsilon}}{\sqrt{t\pi}}e^{\frac{-z^2}{4\varepsilon
t}}
\times\{-\partial_{z}u_{3}^{0}(0,x,y,0)\cos(\alpha t)+\\ \nonumber & &+(\partial_{x}u_{2}^{0}-\partial_{y}u_{1}^{0})(0,x,y,0)\sin(\alpha t)\}\\ \nonumber
& &-\int_{\frac{1}{\sqrt{t}}}^{\infty}\frac{\sqrt{\varepsilon}}{\sqrt{\pi}}\frac{1}{s^2}
e^{\frac{-z^{2}s^{2}}{4\varepsilon}}
\times\{-2\partial_{tz}^{2}u_{3}^{0}(t-\frac{1}{s^2},x,y,0)\cos(\frac{\alpha}{s^{2}})+\\ \nonumber & &+2(\partial_{tx}^{2}u_{2}^{0}-\partial_{ty}^{2}u_{1}^{0})
(t-\frac{1}{s^2},x,y,0)\sin (\frac{\alpha}{s^2})ds\}+
\\ \nonumber & &-
\frac{\sqrt{\varepsilon}}{\sqrt{t\pi}}e^{\frac{-1}{4\varepsilon
t}}
\times\{-\partial_{z}u^{0}_{3}(0,x,y,0)\cos(\alpha t)+\\ \nonumber & & +(\partial_{x}u_{2}^{0}-\partial_{y}u_{1}^{0})(0,x,y,0)\sin(\alpha t)\}-\\ \nonumber
& &-\int_{\frac{1}{\sqrt{t}}}^{\infty}\frac{\sqrt{\varepsilon}}{\sqrt{\pi}}\frac{1}{s^2}
e^{\frac{-s^{2}}{4\varepsilon}}
 \times\{-2\partial_{tz}^{2}u_{3}^{0}(t-\frac{1}{s^2},x,y,0)\cos(\frac{\alpha}{s^{2}})+\\ \nonumber & &+2(\partial_{tx}^{2}u_{2}^{0}-\partial_{ty}^{2}u_{1}^{0})
(t-\frac{1}{s^2},x,y,0)\sin (\frac{\alpha}{s^2})\}ds.
\end{eqnarray}

We denote by $I_{1}+\cdots+I_{4}$ the sum of the terms appearing in the RHS of $(\ref{new1})$.\\ First, we estimate the $L^2$-norm of $I_{1}$, we get:
\begin{eqnarray}\label{estI1}
\|I_{1}\|_{L^{2}(\Omega)}^{2}&\leq&k\varepsilon\int_{0}^{1}e^{\frac{-z^{2}}{2\varepsilon t}}dz\\
&\leq& k\varepsilon\int_{0}^{1}e^{\frac{-cz}{\sqrt{\varepsilon t}}}dz,\hspace{0,5 cm}c>0\nonumber\\
&\leq&k \varepsilon^{3/2}.\nonumber
\end{eqnarray}
Second, we estimate the $L^2$- norm of $I_{2}$, we obtain:
\begin{eqnarray}\label{estI2}
\|I_{2}\|_{L^{2}(\Omega)}^{2}&\leq& k\int_{0}^{1}(\int_{\frac{1}{\sqrt{t}}}^{\infty}\frac{\sqrt{\varepsilon}}{\sqrt{4\pi} s^{2}}e^{\frac{-z^{2}s^2}{4\varepsilon }}ds)^2dz\\
& \leq& (\text{using Cauchy Schwartz inequality})\nonumber\\
&\leq& k\varepsilon \int_{0}^{1}\int_{\frac{1}{\sqrt{t}}}^{\infty}\frac{1}{s^2}ds\int_{\frac{1}{\sqrt{t}}}^{\infty}\frac{1}{s^2}e^{\frac{-z^{2}s^{2}}{2\varepsilon}}dsdz
 \nonumber\\
&\leq& k \varepsilon\int_{0}^{1}\int_{\frac{1}{\sqrt{t}}}^{\infty}\frac{1}{s^2}e^{\frac{-z^{2}s^{2}}{2\varepsilon}}dsdz \nonumber\\
&\leq&k\varepsilon\int_{\frac{1}{\sqrt{t}}}^{\infty}\frac{1}{s^2}\int_{0}^{1}e^{\frac{-csz}{\sqrt{\varepsilon}}}dzds\hspace{0,5cm} c>0,\nonumber\\
&\leq&k\varepsilon^{3/2}. \nonumber
\end{eqnarray}
Finally, combining $(\ref{estI1})$, $(\ref{estI2})$ and the fact that $I_{3}$ and $I_{4}$ are e.s.t. (where \textit{e.s.t.} stands for quantities which are exponentially small terms in all $H^{m}((0,T)\times \Omega),~m\geq 0$), we conclude that:
\begin{equation}\label{estimation}
\left\|\frac{\partial\overline{\varphi}_{3} ^{0,\varepsilon}}{\partial
t}\right\|_{L^{2}(\Omega)}\leq k \varepsilon^{3/4}.
\end{equation}

We will estimate in the following the $z$- derivative of $\overline{\varphi}_{3}
^{0,\varepsilon}$ appearing second term in the RHS of $(\ref{fluide44})_{1}$. Hence by differentiating $\overline{\varphi}_{3}
^{0,\varepsilon}$ with respect to the normal variable $z$, we obtain:
\begin{eqnarray}\label{new2}
\varepsilon\frac{\partial^{2}\overline{\varphi}_{3}^{0,\varepsilon}}{\partial
z^{2}}\!\!\!\!\!&=&\!\!\!\!\!-\!\!\!\int_{0}^{t}\!\!\!\frac{\sqrt{\varepsilon}}{\sqrt{4\pi(t-\tau)}}
\frac{1}{2(t-\tau)}e^{\frac{-z^2}{4\varepsilon (t-\tau)}}
\times\!\!\{2\partial_{z}u^{0}_{3}(\tau,x,y,0)\cos(\alpha(\tau-t))+\\ & &\!\!\!+2(\partial_{x}u_{2}^{0}-\partial_{y}u_{1}^{0})(\tau,x,y,0)\sin(\alpha(\tau-t))\}
d\tau -\nonumber\\
\!\!\!& &-\!\!\!\int_{0}^{t}\!\!\!\!\frac{\sqrt{\varepsilon}}{\sqrt{4\pi(t-\tau)}}
\frac{z^{2}}{4(t-\tau)^{2}}e^{\frac{-z^2}{4\varepsilon
(t-\tau)}}\times\!\!\{2\partial_{z}u^{0}_{3}(\tau,x,y,0)\cos(\alpha(\tau-t))+\nonumber\\
\!\!\!& &\!\!\!+2(\partial_{x}u_{2}^{0}-\partial_{y}u_{1}^{0})(\tau,x,y,0)
\sin(\alpha(\tau-t))
\}d\tau.\nonumber
\end{eqnarray}
We denote by $J_{1}+J_{2}$ the sum of the terms in the RHS of $(\ref{new2})$.
Multiplying $J_{1}$ by $z$ and
setting the change of variables $s=\frac{z}{\sqrt{2\varepsilon(t-\tau)}}$, we find:
\begin{eqnarray*}
z J_{1}&=&\frac{1}{\sqrt{2\pi}}\int_{\frac{z}{\sqrt{2\varepsilon
t}}}^{\infty}\varepsilon
e^{\frac{-s^2}{2}}\times \{-2\partial_{z}u^{0}_{3}(t-\frac{z^2}{2\varepsilon s^2},x,y,0)\cos\left(\frac{\alpha z^2}{2\varepsilon s^2}\right)\\
& &+2(\partial_{x}u_{2}^{0}-\partial_{y}u_{1}^{0})(t-\frac{z^2}{2\varepsilon s^2},x,y,0)
\sin\left(\frac{\alpha z^2}{2\varepsilon s^2}\right)\}ds.
\end{eqnarray*}
Then, we have
\begin{eqnarray*}
|zJ_{1}|
&\leq&k\varepsilon\int_{\frac{z}{\sqrt{2\varepsilon t}}}^{\infty}e^{\frac{-s^2}{2}}ds\\
&\leq&k \varepsilon\int_{\frac{z}{\sqrt{2\varepsilon t}}}^{\infty}e^{-cs}ds,\hspace{0,5cm} c>0\\
&\leq&k \varepsilon e^{\frac{-cz}{\sqrt{2\varepsilon t}}},
\end{eqnarray*}
and we get the $L^{2}$- norm of the term $zJ_{1}$:
\begin{eqnarray*}
\|zJ_{1}\|_{L^{2}(\Omega)}^{2}
&\leq& k \varepsilon^{2}\int_{0}^{1} e^{\frac{-cz}{\sqrt{\varepsilon t}}} dz\\
&\leq& k \varepsilon^{5/2}.
\end{eqnarray*}
Hence, we obtain
\begin{equation}\label{new4}
\|zJ_{1}\|_{L^{2}(\Omega)}\leq k \varepsilon^{5/4}.
\end{equation}
Identically we multiply $J_{2}$ by $z$
and  apply the same change of variables $s=\frac{z}{\sqrt{\varepsilon(t-\tau)}}$, we get:
\begin{eqnarray*}
z J_{2}&=&\frac{1}{\sqrt{\pi}}\int_{\frac{z}{\sqrt{\varepsilon
t}}}^{\infty}2\varepsilon s^{2}
e^{\frac{-s^2}{2}}\times\{-\partial_{z}u^{0}_{3}(t-\frac{z^2}{2\varepsilon s^2},x,y,0)\cos\left(\frac{\alpha z^2}{2\varepsilon s^2}\right)\\ & &+(\partial_{x}u_{2}^{0}-\partial_{y}u_{1}^{0})(t-\frac{z^2}{2\varepsilon s^2},x,y,0)\sin\left(\frac{\alpha z^2}{2\varepsilon s^2})\right)
\}ds.
\end{eqnarray*}
Then, we have
\begin{eqnarray*}
|z J_{2}|&\leq& \int_{\frac{z}{2\sqrt{\varepsilon t}}}^{\infty}k \varepsilon s^{2}e^{\frac{-s^2}{2}}ds\\
&\leq&k\varepsilon (\frac{z}{2\sqrt{\varepsilon t}}e^{\frac{-z^2}{8\varepsilon t}}+\int_{\frac{z}{2\sqrt{\varepsilon t}}}^{\infty}e^{-s^2/2})ds\\
&\leq&k \sqrt{\varepsilon} z e^{\frac{-z^2}{8\varepsilon
t}}+2\varepsilon e^{\frac{-cz}{2\sqrt{\varepsilon t}}},
\end{eqnarray*}
and we obtain the $L^{2}$- norm of $zJ_{2}$,
\begin{eqnarray*}
\|z J_{2}\|_{L^{2}(\Omega)}^{2}&\leq& k\int_{0}^{1}\varepsilon z^{2} e^{\frac{-z^2}{4\varepsilon t}}dz+k \varepsilon^{2}\int_{0}^{1} e^{\frac{-cz}{2\sqrt{\varepsilon t}}}dz\\
&\leq&k\varepsilon\int_{0}^{1}z^{2} e^{\frac{-cz}{\sqrt{2\varepsilon t}}}dz+k \varepsilon^{2}\int_{0}^{1} e^{\frac{-cz}{2\sqrt{\varepsilon t}}}dz\\
&\leq& k \varepsilon^{5/2}.
\end{eqnarray*}
Hence, we infer that
\begin{equation}\label{new3}
\|z J_{2}\|_{L^{2}(\Omega)}\leq k \varepsilon^{5/4}.
\end{equation}
Finally, combining $(\ref{new4})$ and $(\ref{new3})$, we deduce the
following estimate:
\begin{equation}\label{estimation1}
\|z\varepsilon\frac{\partial^{2}\overline{\varphi}_{3}^{0,\varepsilon}}{\partial
z^{2}}\|_{L^{2}(\Omega)}\leq k \varepsilon^{5/4}.
\end{equation}
\section{Convergence result}
In this section we prove the main theoretical result of this article.
\begin{theorem}\label{convergence}
The solution $\boldsymbol{u}^{\varepsilon}$ of $(\ref{fluide1})$, with $\boldsymbol{u}_{0}$ and $\boldsymbol{f}$
supposed to be sufficiently smooth, satisfies the following
estimates:
\begin{equation}\label{t1}
\|\boldsymbol{u}^{\varepsilon}-\boldsymbol{u}^{0}-\overline{\boldsymbol{\varphi}}^{0,\varepsilon}
-\widetilde{\boldsymbol{\varphi}}^{1,\varepsilon}\|_{L^{\infty}(0,T,\textbf{L}^{2}(\Omega))}\leq
k \varepsilon^{3/4},
\end{equation}
\begin{equation}\label{t2}
\|\boldsymbol{u}^{\varepsilon}-\boldsymbol{u}^{0}-\overline{\boldsymbol{\varphi}}^{0,\varepsilon}-
\widetilde{\boldsymbol{\varphi}}^{1,\varepsilon}\|_{L^{\infty}(0,T,\textbf{H}^{1}(\Omega))}\leq
k \varepsilon^{1/4},
\end{equation}
where $k$ is a positive constant depending on the data but not
$\varepsilon$, and  $\boldsymbol{u}^{0}$ and $\boldsymbol{\varphi}^{\varepsilon}$ are
defined respectively by $(\ref{eq lim 00})$ and $(\ref{fluide
3})$. Here we denoted by $\textbf{L}^{2}(\Omega)=(L^{2}(\Omega))^{3}$ and
$\textbf{H}^{1}(\Omega)=(H^{1}(\Omega))^{3}.$
\end{theorem}
\begin{preuve}
\normalfont First we observe that the corrector
$\boldsymbol{\varphi}^{\varepsilon}$ does not satisfy the desired boundary
conditions as given by $(\ref{fluide 3})_{3}$, this is due to the choice
of a corrector in a simpler form. To overcome this difficulty we
introduce additional (small) correctors $\overline{\boldsymbol{\theta}}^{\varepsilon}$ and $\widetilde{\boldsymbol{\theta}}^{\varepsilon}$ as follows:
\begin{equation}\label{fluidestokes}
\left\lbrace
\begin{array}{lcl}
-\varepsilon\Delta\overline{\boldsymbol{\theta}}^{\varepsilon}+\nabla\Pi^{\varepsilon}=0, \hspace{0,5cm}\mbox{in}\hspace{0,2cm}\Omega_{\infty}\times(0,T),\\
\textnormal{div}\overline{\boldsymbol{\theta}}^{\varepsilon}=0,\\
\overline{\boldsymbol{\theta}}^{\varepsilon}|_{z=0}=(0,0,-\overline{\varphi}_{3}^{0,\varepsilon}|_{z=0}),\\
\overline{\boldsymbol{\theta}}^{\varepsilon}|_{z=1}=(-\overline{\varphi}_{1}^{0,\varepsilon}|_{z=1},-\overline{\varphi}_{2}^{0,\varepsilon}|_{z=1},0),\\
\end{array}\right.
\end{equation}
and
\begin{equation}\label{fluidestokes1}
\left\lbrace
\begin{array}{lcl}
-\varepsilon\Delta\widetilde{\boldsymbol{\theta}}^{\varepsilon}+\nabla Q^{\varepsilon}=0, \hspace{0,5cm}\mbox{in}\hspace{0,2cm}\Omega_{\infty}\times(0,T),\\
\textnormal{div}\widetilde{\boldsymbol{\theta}}^{\varepsilon}=0,\\
\widetilde{\boldsymbol{\theta}}^{\varepsilon}|_{z=1}=(0,0,-\widetilde{\varphi}_{3}^{0,\varepsilon}|_{z=1}),\\
\widetilde{\boldsymbol{\theta}}^{\varepsilon}|_{z=0}=(-\widetilde{\varphi}_{1}^{0,\varepsilon}|_{z=0},-\widetilde{\varphi}_{2}^{0,\varepsilon}|_{z=0},0).\\
\end{array}\right.
\end{equation}
To estimate the $L^2$- norm of the additional correctors, we set
$\overline{\boldsymbol{\theta}}^{\varepsilon}=\sqrt{\varepsilon}\widetilde{\overline{\boldsymbol{\theta}}}^{\varepsilon},
\Pi^{\varepsilon}=\varepsilon^{3/2}\widetilde{\Pi^{\varepsilon}},$ hence $\widetilde{\overline{\boldsymbol{\theta}}}^{\varepsilon}$ satisfies the following system:
\begin{equation}\label{fluide}
\left\lbrace
\begin{array}{lcl}
-\Delta\widetilde{\overline{\boldsymbol{\theta}}}^{\varepsilon}+\nabla\widetilde{\Pi^{\varepsilon}}=0, \hspace{0,5cm}\mbox{in}\hspace{0,2cm}\Omega_{\infty}\times(0,T)\\
\textnormal{div}\widetilde{\overline{\boldsymbol{\theta}}}^{\varepsilon}=0,\\
\overline{\boldsymbol{\theta}}^{\varepsilon}|_{z=0}=(0,0,-\frac{\overline{\varphi}_{3}^{0,\varepsilon}}{\sqrt{\varepsilon}}|_{z=0}),\\
\overline{\boldsymbol{\theta}}^{\varepsilon}|_{z=1}=(-\frac{\overline{\varphi}_{1}^{0,\varepsilon}}{\sqrt{\varepsilon}}|_{z=1},
-\frac{\overline{\varphi}_{2}^{0,\varepsilon}}{\sqrt{\varepsilon}}|_{z=1},0).\\
\end{array}\right.
\end{equation}
Then we deduce from the direct estimates of the Stokes problem
(see \cite{1} ) that:
\begin{eqnarray*}
\|\widetilde{\overline{\boldsymbol{\theta}}}^{\varepsilon}\|_{L^{2}(\Omega)}&\leq&k\|\frac{\overline{\varphi}_{3}^{0,\varepsilon}}
{\sqrt{\varepsilon}}|_{z=0}
\|_{H^{-1/2}(\Gamma)}
+k\|\frac{\overline{\varphi}_{1}^{0,\varepsilon}}{\sqrt{\varepsilon}}|_{z=1}\|_{H^{-1/2}(\Gamma)}+
k\|\frac{\overline{\varphi}_{2}^{0,\varepsilon}}{\sqrt{\varepsilon}}|_{z=1}\|_{H^{-1/2}(\Gamma)}\\
&\leq&k\|\frac{\overline{\varphi}_{3}^{0,\varepsilon}}{\sqrt{\varepsilon}}\|_{L^{2}(\Omega)}+e.s.t.
\end{eqnarray*}
Now we will estimate the $ L^{2}$- norm of
$\dfrac{\overline{\varphi}_{3}^{0,\varepsilon}}{\sqrt{\varepsilon}}$,
hence we have:

\begin{eqnarray*}
|\frac{\overline{\varphi}_{3}^{0,\varepsilon}}{\sqrt{\varepsilon}}|^{2}&\leq&k(\int_{\frac{1}{\sqrt{t}}}^{\infty}\frac{1}{4\sqrt{\pi}s^{2}}
e^{\frac{-z^{2}s^{2}}{4\varepsilon}}ds)^{2}\\
&\leq & (\text{Using the Cauchy-Schwartz inequality})\\
&\leq&k\int_{\frac{1}{\sqrt{t}}}^{\infty}\frac{1}{s^{2}}ds\int_{\frac{1}{\sqrt{t}}}^{\infty}\frac{1}{s^{2}}e^{\frac{-z^{2}s^{2}}{2\varepsilon}}ds\\
&\leq&k\int_{\frac{1}{\sqrt{t}}}^{\infty}\frac{1}{s^{2}}e^{\frac{-z^{2}s^{2}}{2\varepsilon}}ds.
\end{eqnarray*}
Therefore, we have
\begin{eqnarray*}
\|\frac{\overline{\varphi}_{3}^{0,\varepsilon}}{\sqrt{\varepsilon}}\|_{L^{2}(\Omega)}^{2}&\leq&k\int_{0}^{1}
\int_{\frac{1}{\sqrt{t}}}^{\infty}\frac{1}{s^2}e^{\frac{-z^{2}s^{2}}{2\varepsilon}}dsdz\\
&\leq&k\int_{\frac{1}{\sqrt{t}}}^{\infty}\frac{1}{s^{2}}\int_{0}^{1}e^{\frac{-z^{2}s^{2}}{2\varepsilon}}dzds\\
&\leq&k\int_{\frac{1}{\sqrt{t}}}^{\infty}\frac{1}{s^{2}}\int_{0}^{1}e^{\frac{-czs}{\sqrt{\varepsilon}}}dzds,\hspace{0,5cm}c>0\\
&\leq&k\sqrt{\varepsilon},
\end{eqnarray*}
Hence, we infer that
\begin{eqnarray*}
\|\widetilde{\overline{\boldsymbol{\theta}}}^{\varepsilon}\|_{L^{2}(\Omega)}&\leq&k
\varepsilon^{1/4}.
\end{eqnarray*}
Finally, we get
\begin{equation}\label{estimation2}
\|\overline{\boldsymbol{\theta}}^{\varepsilon}\|_{L^{2}(\Omega)}\leq k
\varepsilon^{3/4}.
\end{equation}
In the following we will estimate the $L^2(\Omega)$ norm of the gradient of
$\widetilde{\overline{\boldsymbol{\theta}}}^{\varepsilon}$, hence we
find:
\begin{eqnarray*}
\|\nabla\widetilde{\overline{\boldsymbol{\theta}}}^{\varepsilon}\|_{L^{2}(\Omega)}&\leq&k\|
\frac{\overline{\varphi}_{3}^{0,\varepsilon}}{\sqrt{\varepsilon}}|_{z=0}\|_{H^{1/2}(\Gamma)}
+k\|\frac{\overline{\varphi}_{1}^{0,\varepsilon}}{\sqrt{\varepsilon}}|_{z=1}\|_{H^{1/2}(\Gamma)}+
k\|\frac{\overline{\varphi}_{2}^{0,\varepsilon}}{\sqrt{\varepsilon}}|_{z=1}\|_{H^{1/2}(\Gamma)}\\
&\leq&k\|\frac{\overline{\varphi}_{3}^{0,\varepsilon}}{\sqrt{\varepsilon}}\|_{H^{1}(\Omega)}+e.s.t\\
&\leq&k \varepsilon^{-1/4}.
\end{eqnarray*}
Thus we deduce that:
\begin{equation}\label{estimation2'}
\|\nabla\overline{\boldsymbol{\theta}}^{\varepsilon}\|_{L^{2}(\Omega)}\leq k
\varepsilon^{1/4}.
\end{equation}
 We notice that the estimate $(\ref{estimation2})$ also holds for
 the time derivative of $\overline{{\theta}}^{\varepsilon}$, i.e.,
\begin{equation}\label{estimation3}
\|\frac{\partial\overline{\boldsymbol{\theta}}^{\varepsilon}}{\partial
t}\|_{L^{2}(\Omega)}\leq k \varepsilon^{3/4}.
\end{equation}
We now define
$\boldsymbol{w}^{\varepsilon}=\boldsymbol{u}^{\varepsilon}-\boldsymbol{u}^{0}-\overline{\boldsymbol{\varphi}}^{0,\varepsilon}
-\widetilde{\boldsymbol{\varphi}}^{1,\varepsilon}-\overline{\boldsymbol{\theta}}^{\varepsilon}-\widetilde
{\boldsymbol{\theta}}^{\varepsilon}$, and according to
$(\ref{fluide1}), (\ref{eq lim 00}), (\ref{fluide44}), (\ref{fluidestokes})$ and $(\ref{fluidestokes1}),~\boldsymbol{w}^{\varepsilon}$ verifies:\\
\begin{equation}\label{fluide fin}
\left\lbrace
\begin{array}{lcl}
\dfrac{\partial\boldsymbol{w}^{\varepsilon}}{\partial t}-\varepsilon\Delta
\boldsymbol{w}^{\varepsilon}+\boldsymbol{\omega}\times \boldsymbol{w}^{\varepsilon}
+\nabla(p^{\varepsilon}-p^{0}-\Pi^{\varepsilon}-Q^{\varepsilon})=\varepsilon\dfrac{\partial^{2}\overline{\boldsymbol{\varphi}}^{0,\varepsilon}}{\partial
x^{2}}+\varepsilon\dfrac{\partial^{2}\widetilde{\boldsymbol{\varphi}}^{1,\varepsilon}}{\partial
x^{2}}\\+\varepsilon\dfrac{\partial^{2}\overline{\boldsymbol{\varphi}}^{0,\varepsilon}}{\partial
y^{2}}+\varepsilon\dfrac{\partial^{2}\widetilde{\boldsymbol{\varphi}}^{1,\varepsilon}}{\partial
y^{2}}+\varepsilon\Delta
\boldsymbol{u}^{0}-\boldsymbol{\omega}\times\overline{\boldsymbol{\theta}}^{\varepsilon}-
\boldsymbol{\omega}\times\widetilde{\boldsymbol{\theta}}^{\varepsilon}-\dfrac{\partial\overline{\boldsymbol{\theta}}^{\varepsilon}}{\partial
t} -\dfrac{\partial\widetilde{\boldsymbol{\theta}}^{\varepsilon}}{\partial
t}\\+(0,0,\dfrac{\partial
\overline{\varphi}_{3}^{0,\varepsilon}}{\partial
t})+(0,0,\dfrac{\partial
\widetilde{\varphi}_{3}^{1,\varepsilon}}{\partial
t})+(0,0,\varepsilon\dfrac{\partial^{2}\overline{\varphi}_{3}^{0,\varepsilon}}{\partial z^{2}})+
(0,0,\varepsilon\dfrac{\partial^{2}\widetilde{\varphi}_{3}^{1,\varepsilon}}{\partial z^{2}}), \hspace{0,5cm}\mbox{in}\hspace{0,2cm}\Omega_{\infty}\times(0,T),\vspace{0,2cm}\\
\textnormal{div}~\boldsymbol{w}^{\varepsilon}=0,\mbox{in}\hspace{0,2cm}\Omega_{\infty}\times(0,T),\vspace{0,2cm}\\
\boldsymbol{w}^{\varepsilon}=0, \hspace{0.5cm} \text{at}~~z=0,1,\\
\boldsymbol{w}^{\varepsilon}~\text{is}~2\pi\text{-periodic in the}~ x~
\text{and}~y~\text{directions},\\
 \boldsymbol{w}^{\varepsilon}|_{t=0}=0.
\end{array}\right.
\end{equation}
We multiply $(\ref{fluide fin})_{1}$ by $\boldsymbol{w}^{\varepsilon}$, integrate over $\Omega$, and apply the Cauchy-Shwarz inequality, we obtain:\\

\begin{eqnarray*}
\dfrac{1}{2}\dfrac{d\|\boldsymbol{w}^{\varepsilon}\|^{2}}{dt}\!+\!\varepsilon\|\nabla
\boldsymbol{w}^{\varepsilon}\|^{2}\!\!\!&\leq&\!\!\!\!\varepsilon\|\frac{\partial^{2}\overline{\boldsymbol{\varphi}}^{0,\varepsilon}}{\partial
x^{2}}\|\|\boldsymbol{w}^{\varepsilon}\|\!+\!\varepsilon\|\frac{\partial^{2}\widetilde{\boldsymbol{\varphi}}^{1,\varepsilon}}{\partial
x^{2}}\|\|\boldsymbol{w}^{\varepsilon}\|\!+\!\varepsilon\|\frac{\partial^{2}\overline{\boldsymbol{\varphi}}^{0,\varepsilon}}{\partial
y^{2}}\|\|\boldsymbol{w}^{\varepsilon}\|+\\ \nonumber& &\!\!\!\!+\varepsilon\|\frac{\partial^{2}\widetilde{\boldsymbol{\varphi}}^{0,\varepsilon}}{\partial
y^{2}}\|\|\boldsymbol{w}^{\varepsilon}\|\!+\!\varepsilon\|\Delta
u^{0}\|\|\boldsymbol{w}^{\varepsilon}\|\!+\!\|\overline{\boldsymbol{\theta}}^{\varepsilon}\|\|\boldsymbol{w}^{\varepsilon}\| \!+\!\|\widetilde{\boldsymbol{\theta}}^{\varepsilon}\|\|\boldsymbol{w}^{\varepsilon}\|+\\ \nonumber& &\!\!\!\!+\|\frac{\partial\overline{\boldsymbol{\theta}}
^{\varepsilon}}{\partial
t}\|\|\boldsymbol{w}^{\varepsilon}\|\!+\!\|\frac{\partial\widetilde{\boldsymbol{\theta}}
^{\varepsilon}}{\partial t}\|\|\boldsymbol{w}^{\varepsilon}\|\!+\!
\|\varepsilon z\frac{\partial^{2}\overline{\varphi}_{3}^{0,\varepsilon}}{\partial^{2}z}\|\|\nabla w^{\varepsilon}\|
+\\ \nonumber & &\!\!\!\!+\|\varepsilon z\frac{\partial^{2}\widetilde{\varphi}_{3}^{0,\varepsilon}}{\partial^{2}z}\|\|\nabla \boldsymbol{w}^{\varepsilon}\|
\!+\!\|\frac{\partial \overline{\varphi}_{3}^{1,\varepsilon}}{\partial t}\|\|\boldsymbol{w}^{\varepsilon}\|
\!+\!\|\frac{\partial\widetilde{\varphi}_{3}^{1,\varepsilon}}{\partial t}\|\|\boldsymbol{w}^{\varepsilon}\|.
\end{eqnarray*}
Hence according to $(\ref{estimation}), (\ref{estimation1}),
(\ref{estimation2})$ and $(\ref{estimation3})$, we have:
\begin{eqnarray*}
\dfrac{1}{2}\dfrac{d\|\boldsymbol{w}^{\varepsilon}\|}{dt}+\varepsilon\|\nabla \boldsymbol{w}^{\varepsilon}\|^{2}&\leq&\frac{1}{2}\|\boldsymbol{w}^{\varepsilon}\|^{2}
+k \varepsilon^{3/2}+k \varepsilon^{3/4}\varepsilon^{1/2}\frac{\|\nabla\boldsymbol{w}^{\varepsilon}\|}{2}+k \varepsilon^{3/4}\varepsilon^{1/2}\frac{\|\nabla \boldsymbol{w}^{\varepsilon}\|}{2}\\
&\leq&\frac{1}{2}\|\boldsymbol{w}^{\varepsilon}\|^{2}+k
\varepsilon^{3/2}+\frac{\varepsilon}{2}\|\nabla
\boldsymbol{w}^{\varepsilon}\|^{2}.
\end{eqnarray*}
In conclusion, we have
$$\dfrac{d\|\boldsymbol{w}^{\varepsilon}\|^{2}}{dt}+\varepsilon
\|\nabla \boldsymbol{w}^{\varepsilon}\|^{2}\leq \|\boldsymbol{w}^{\varepsilon}\|^{2}+k
\varepsilon^{3/2}.$$
Using the Gronwall inequality, we obtain\\
$$\|\boldsymbol{w}^{\varepsilon}\|_{L^{\infty}(0,T;\boldsymbol{L^{2}}(\Omega))}\leq k \varepsilon^{3/4}~\text{and}~\|\nabla \boldsymbol{w}^{\varepsilon}\|_{L^{\infty}(0,T;\boldsymbol{L^{2}}(\Omega))}\leq k \varepsilon^{1/4}.$$
Hence, according to $(\ref{estimation2}), (\ref{estimation2'}),$
and the triangular inequality, we deduce $(\ref{t1})$ and $
(\ref{t2})$. This concludes the proof of Theorem \ref{convergence}.
\end{preuve}
\section{A collocated finite volume scheme with a splitting method
for the time discretization} We follow here the notations of \cite{5n} that we recall in this section for the reader convenience. In the following, we uniformly discretize the domain $\Omega$
by using cube finite volumes of dimensions $\Delta x
\Delta y \Delta z$:
$$K_{i,j,k}=[x_{i-\frac{1}{2}},x_{i+\frac{1}{2}}]\times[y_{j-\frac{1}{2}},y_{j+\frac{1}{2}}]\times
[z_{k-\frac{1}{2}},z_{k+\frac{1}{2}}],$$ where:
$$x_{i+\frac{1}{2}}=i\Delta x,\quad y_{j+\frac{1}{2}}=j\Delta y,\quad z_{k+\frac{1}{2}}=k\Delta
z,$$ $$\forall~i=0,\ldots,M, \forall~j=0,\ldots N,
\forall~k=0,\ldots L.$$
The edges of the control volumes are defined by:
$$\Gamma_{i+1/2,j,k}=\{(x,y,z);x=x_{i+1/2},y\in [y_{j-1/2},y_{j+1/2}], z\in [z_{k-1/2}, z_{k+1/2}]\},$$
$$\Gamma_{i,j+1/2,k}=\{(x,y,z);x\in [x_{i-1/2},x_{i+1/2}], y=y_{j+1/2}, z\in [z_{k-1/2}, z_{k+1/2}]\},$$
$$\Gamma_{i,j,k+1/2}=\{(x,y,z);x\in [x_{i-1/2},x_{i+1/2}], y\in [y_{j-1/2},y_{j+1/2}], z=z_{k+1/2}]\},$$
$$\forall~i=0,\ldots,M, \forall~j=0,\ldots N,
\forall~k=0,\ldots L.$$
The velocity and the pressure are approximated in the center of the cells as follows:
$$\boldsymbol{u}_{i,j,k}(t)\simeq\dfrac{1}{\Delta x\Delta y \Delta
z}\int_{x_{i-\frac{1}{2}}}^{x_{i+\frac{1}{2}}}\int_{y_{j-\frac{1}{2}}}^{y_{j+\frac{1}{2}}}\int_{z_{k-\frac{1}{2}}}^{z_{k+\frac{1}{2}}}\boldsymbol{u}(x,y,z,t)dxdydz,$$
$$p_{i,j,k}(t)\simeq\dfrac{1}{\Delta x\Delta y \Delta
z}\int_{x_{i-\frac{1}{2}}}^{x_{i+\frac{1}{2}}}\int_{y_{j-\frac{1}{2}}}^{y_{j+\frac{1}{2}}}\int_{z_{k-\frac{1}{2}}}^{z_{k+\frac{1}{2}}}p(x,y,z,t)dxdydz.$$
We also define the velocity fluxes:
$$F_{u i+\frac{1}{2},j,k}\simeq\frac{1}{\Delta
y\Delta z}\int_{y_{j-\frac{1}{2}}}^{y_{j+\frac{1}{2}}}\int_{z_{k-\frac{1}{2}}}^{z_{k+\frac{1}{2}}}u(x_{i+\frac{1}{2}},y,z,t)dy
dz,$$

$$F_{v i,j+\frac{1}{2},k}\simeq\frac{1}{\Delta
x\Delta z}\int_{x_{i-\frac{1}{2}}}^{x_{i+\frac{1}{2}}}\int_{z_{k-\frac{1}{2}}}^{z_{k+\frac{1}{2}}}v(x,y_{j+\frac{1}{2}},z,t)dx
dz,$$

$$F_{w i,j,k+\frac{1}{2}}\simeq\frac{1}{\Delta x\Delta
y}\int_{x_{i-\frac{1}{2}}}^{x_{i+\frac{1}{2}}}\int_{y_{j-\frac{1}{2}}}^{y_{j+\frac{1}{2}}}w(x,y,z_{k+\frac{1}{2}},t)dx
dy.$$

\subsection{Time disceretization}
We start by choosing a time disceretization for
$(\ref{fluide1})_{1}$:

\begin{equation}\label{num}
\dfrac{3 \boldsymbol{u}^{n+1}-4 \boldsymbol{u}^{n}+\boldsymbol{u}^{n-1}}{2\Delta t}-\varepsilon\Delta
\boldsymbol{u}^{n+1}+2\boldsymbol{\omega}\times \boldsymbol{u}^{n}-\boldsymbol{\omega}\times
\boldsymbol{u}^{n-1}+2\nabla p^{n}-\nabla p^{n-1}=\boldsymbol{f}^{n+1}.
\end{equation}
Thanks to ($\ref{num}$) we are able to compute the new velocity
$\boldsymbol{u}^{n+1}$.\\
Hence, to obtain the pressure, we take the divergence of
$(\ref{fluide1})_{1}$ and use the incompressibility condition $(\ref{fluide1})_{2}$ we find:
\begin{equation}\label{rr}
\Delta p=\textnormal{div}(\boldsymbol{f}+\varepsilon\Delta \boldsymbol{u}+\boldsymbol{\omega}\times
\boldsymbol{u}).
\end{equation}
Thus we discretize $(\ref{rr})$ as follows:
\begin{equation}\label{numerique}
\Delta p^{n+1}=\textnormal{div}(\boldsymbol{f}^{n+1}+\varepsilon\Delta
\boldsymbol{u}^{n+1}-2\boldsymbol{\omega}\times
\boldsymbol{u}^{n}-\boldsymbol{\omega}\times \boldsymbol{u}^{n-1}).
\end{equation}
By replacing $\Delta$ by $-\nabla\times\nabla\times$ (see\cite{5n}~\text{and}~\cite{7n}), we rewrite $(\ref{numerique})$ as below:
\begin{equation}\label{ch}
\Delta p^{n+1}=\textnormal{div}(\boldsymbol{f}^{n+1}-\varepsilon\nabla\times\nabla\times
\boldsymbol{u}^{n+1}-2\boldsymbol{\omega}\times
\boldsymbol{u}^{n}-\boldsymbol{\omega}\times \boldsymbol{u}^{n-1}).
\end{equation}
Now, by using the relation $ \Delta u=\nabla \textnormal{div} \boldsymbol{u}^{n+1}-\nabla \times \nabla \times \boldsymbol{u}^{n+1}$,
then $(\ref{num})$ becomes:
\begin{eqnarray*}
\boldsymbol{f}^{n+1}-\varepsilon\nabla\times\nabla\times
\boldsymbol{u}^{n+1}-2\omega\times \boldsymbol{u}^{n}-\boldsymbol{\omega}\times
\boldsymbol{u}^{n-1}\\=\dfrac{3 \boldsymbol{u}^{n+1}-4 \boldsymbol{u}^{n}+\boldsymbol{u}^{n-1}}{2\Delta
t}-\varepsilon\nabla \textnormal{div} \boldsymbol{u}^{n+1}+2\nabla p^{n}-\nabla
p^{n-1}.\nonumber
\end{eqnarray*}
Hence, we deduce from $(\ref{ch})$ that
\begin{equation}
\Delta p^{n+1}=\textnormal{div}(\dfrac{3\boldsymbol{u}^{n+1}-4\boldsymbol{u}^{n}+\boldsymbol{u}^{n-1}}{2\Delta
t}-\varepsilon\nabla \textnormal{div} \boldsymbol{u}^{n+1}+2\nabla p^{n}-\nabla
p^{n-1}).
\end{equation}
Thus, we obtain
\begin{equation}
\Delta(p^{n+1}-2p^{n}+p^{n-1}+\varepsilon \textnormal{div}
\boldsymbol{\^u}^{n+1})=\textnormal{div}(\dfrac{3\boldsymbol{u}^{n+1}-4\boldsymbol{u}^{n}+\boldsymbol{u}^{n-1}}{2\Delta t}).
\end{equation}
Then we compute the pressure from
\begin{equation}\label{numerique3}
\left\lbrace
\begin{array}{lcl}
\Delta \psi^{n+1}=\textnormal{div}(\dfrac{3\boldsymbol{u}^{n+1}-4\boldsymbol{u}^{n}+\boldsymbol{u}^{n-1}}{2\Delta t}),\\
\displaystyle\frac{\partial\psi^{n+1}}{\partial n}=0,\\
\end{array}\right.
\end{equation}
and
\begin{equation}
p^{n+1}=\psi^{n+1}+2 p^{n}-p^{n-1}-\varepsilon\textnormal{div} \boldsymbol{u}^{n+1}.
\end{equation}
Concerning the boundary conditions, we have the periodicity in the
$x$ and $y$ directions and the Dirichlet boundary conditions in the $z$
direction for $u^{n+1}$:
$$ \boldsymbol{u}^{n+1}_{0,j,k}=\boldsymbol{u}^{n+1}_{N,j,k},\quad\quad
\boldsymbol{u}^{n+1}_{N+1,j,k}=\boldsymbol{u}^{n+1}_{1,j,k},$$
$$ \boldsymbol{u}^{n+1}_{i,0,k}=\boldsymbol{u}^{n+1}_{i,N,k},\quad\quad
\boldsymbol{u}^{n+1}_{i,N+1,k}=\boldsymbol{u}^{n+1}_{i,1,k},$$
$$\dfrac{\boldsymbol{u}^{n+1}_{i,j,N+1}+\boldsymbol{u}^{n+1}_{i,j,N}}{2}=0,\quad\quad\dfrac{\boldsymbol{u}^{n+1}_{i,j,0}+\boldsymbol{u}^{n+1}_{i,j,1}}{2}=0.$$
The Neumann boundary conditions are imposed for $\psi^{n+1}$ in the $z$ direction
and the periodicity in $x$ and $y$ directions. Thus, we have
$$ \psi^{n+1}_{0,j,k}=\psi^{n+1}_{N,j,k},\quad\quad
\psi^{n+1}_{N+1,j,k}=\psi^{n+1}_{1,j,k},$$
$$ \psi^{n+1}_{i,0,k}=\psi^{n+1}_{i,N,k},\quad\quad
\psi^{n+1}_{i,N+1,k}=\psi^{n+1}_{i,1,k},$$
$$\psi^{n+1}_{i,j,N+1}=\psi^{n+1}_{i,j,N},\quad\quad\psi^{n+1}_{i,j,0}=\psi^{n+1}_{i,j,1}.$$
The periodicity in $x$ and $y$ for the pressure yields:
$$p_{0,j,k}=p_{M,j,k},\quad \quad p_{M+1,j,k}, p_{1,j,k},$$
$$p_{i,0,k}=p_{i,M,k},\quad \quad p_{i,M+1,k}=p_{i,1,k},$$
and for the terms $p_{i,j,0}$ and $p_{i,j,L+1}$ we use the second order compact scheme to compute them:
$$p_{i,j,0}=\frac{5}{2}p_{i,j,1}-2p_{i,j,2}+\frac{1}{2}p_{i,j,3},\quad \quad p_{i,j,L+1}=\frac{5}{2}p_{i,j,L}-2p_{i,j,L-1}+\frac{1}{2}p_{i,j,L-2}.$$

\subsection{Finite volume discretization}
To Compute the velocity $\boldsymbol{u}^{n+1}$, we discretize $(\ref{num})$ and we obtain:
\begin{eqnarray*} \label{tem}&
&\Delta x\Delta y \Delta z\dfrac{3\boldsymbol{u}^{n+1}-4\boldsymbol{u}^{n}+\boldsymbol{u}^{n-1}}{2\Delta
t}-\varepsilon [\Delta x \Delta
y\dfrac{\boldsymbol{u}^{n+1}_{i,j,k+1}-2 \boldsymbol{u}^{n+1}_{i,j,k}+\boldsymbol{u}^{n-1}_{i,j,k-1}}{\Delta z}\\&
&+\Delta y \Delta
z\dfrac{\boldsymbol{u}^{n+1}_{i+1,j,k}-2 \boldsymbol{u}^{n+1}_{i,j,k}+\boldsymbol{u}^{n+1}_{i-1,j,k}}{\Delta x}+\Delta
x \Delta z\dfrac{\boldsymbol{u}^{n+1}_{i,j+1,k}-2 \boldsymbol{u}^{n+1}_{i,j,k}+\boldsymbol{u}^{n+1}_{i,j-1,k}}{\Delta
y}]\\&+&2\left(\begin{array}{c}
\dfrac{\Delta y\Delta z}{2}(p^{n}_{i+1,j,k}-p^{n}_{i-1,j,k}) \\\\
\dfrac{\Delta x\Delta z}{2}(p^{n}_{i,j+1,k}-p^{n}_{i,j-1,k})\\\\
\dfrac{\Delta x\Delta y}{2}(p^{n}_{i,j,k+1}-p^{n}_{i-1,j,k-1}) \\
\end{array}\right)
-\left(\begin{array}{c}
\dfrac{\Delta y\Delta z}{2}(p^{n-1}_{i+1,j,k}-p^{n-1}_{i-1,j,k}) \\\\
\dfrac{\Delta x\Delta z}{2}(p^{n-1}_{i,j+1,k}-p^{n-1}_{i,j-1,k})  \\\\
\dfrac{\Delta x\Delta y}{2}(p^{n-1}_{i,j,k+1}-p^{n-1}_{i-1,j,k-1})\\
  \end{array}\right)
\\&+&\Delta x \Delta y \Delta z (\boldsymbol{\omega}\times (2\boldsymbol{u}^{n}_{i,j,k}-\boldsymbol{u}^{n-1}_{i,j,k}))=\Delta
x\Delta y \Delta z \boldsymbol{f}^{n+1}_{i,j,k}.
\end{eqnarray*}

To compute the pressure we first compute $\psi^{n+1}$:

\begin{eqnarray*}
& &\Delta x \Delta y\dfrac{\psi^{n+1}_{i,j,k+1}-2\psi^{n+1}_{i j k}+\psi^{n+1}_{i
j k-1}}{\Delta z}+\Delta y \Delta
z\dfrac{\psi^{n+1}_{i+1,j,k}-2\psi^{n+1}_{i,j,k}+\psi^{n+1}_{i-1,j,k}}{\Delta
x}\\& &+\Delta x \Delta
z\dfrac{\psi^{n+1}_{i,j+1,k}-2\psi^{n+1}_{i,j,k}+\psi^{n+1}_{i,j-1,k}}{\Delta
y}=\frac{1}{2\Delta t}[\Delta y \Delta z[(3F^{n+1}_{u
i+\frac{1}{2} j k}-4 F^{n}_{u i+\frac{1}{2} j k}+F^{n-1}_{u
i+\frac{1}{2} j k})\\& &-(3F^{n+1}_{u i-\frac{1}{2} j k}-4
F^{n}_{u i-\frac{1}{2} j k}+F^{n-1}_{u i-\frac{1}{2} j k})]+\Delta
x \Delta z[(3F^{n+1}_{v i j+\frac{1}{2} k}-4 F^{n}_{v i
j+\frac{1}{2} k}+F^{n-1}_{v i j+\frac{1}{2} k})\\& &-(3F^{n+1}_{v
i j-\frac{1}{2} k}-4 F^{n}_{v i j-\frac{1}{2} k}+F^{n-1}_{v i
j-\frac{1}{2} k})]+\Delta x \Delta y[(3F^{n+1}_{w i j
k+\frac{1}{2}}-4 F^{n}_{w i j k+\frac{1}{2}}+F^{n-1}_{w i j
k+\frac{1}{2}})\\& &-(3F^{n+1}_{w i j k-\frac{1}{2} k}-4 F^{n}_{w
i j k-\frac{1}{2}}+F^{n-1}_{w i j k-\frac{1}{2}})]].
\end{eqnarray*}

Then, we easily obtain the pressure:\\
\begin{eqnarray*}
p^{n+1}_{i,j,k}&=&\psi^{n+1}_{i,j,k}+2p^{n}_{i,j,k}-p^{n-1}_{i,j,k}-\dfrac{\varepsilon}{\Delta
x \Delta y \Delta z}[\Delta y \Delta z(F^{n+1}_{u i+\frac{1}{2} j
k}-F^{n+1}_{u i-\frac{1}{2} j k})\\& &+\Delta x \Delta
z(F^{n+1}_{v i j+\frac{1}{2} k}-F^{n+1}_{v i j-\frac{1}{2}
k})+\Delta x \Delta y(F^{n+1}_{w i j k+\frac{1}{2}}-F^{n+1}_{w i j
k-\frac{1}{2}})].
\end{eqnarray*}
\subsection{Computation of the fluxes}
We recall here that the simplest method to compute the fluxes (linear interpolation)
does not work when the viscosity $\varepsilon$ is small. Hence the authors in \cite{5n} considered a modified interpolation method for the fluxes in two dimensional case. Now, since we aim here to study the boundary layers at small viscosity, we need, on the one hand, to adapt the discretization in \cite{5n} to the $3D$ dimensional case and, on the other hand, to introduce the correctors in the finite volume discretization basis that is the NFVM. Thus we first start by introducing the $3D$ fluxes inherited from \cite{5n}:
\begin{eqnarray*}
F_{u_{ i+\frac{1}{2},j,k}^{n+1}}&=&\frac{u_{i+1,j,k}^{n+1}+u_{i,j,k}^{n+1}}{2}+\theta\frac{\Delta y\Delta z}{4a}(p_{i+2,j,k}^{n}-2p_{i+1,j,k}^{n}+p_{i,j,k}^{n})\\& &-
\theta\frac{\Delta y\Delta z}{4a}(p_{i+1,j,k}^{n}-2p_{i,j,k}^{n}+p_{i-1,j,k}^{n}),
\end{eqnarray*}
\begin{eqnarray*}
F_{v _{i+\frac{1}{2},j,k}^{n+1}}&=&\frac{v_{i,j+1,k}^{n+1}+v_{i,j,k}^{n+1}}{2}+\theta\frac{\Delta x\Delta z}{4a}(p_{i,j+2,k}^{n}-2p_{i,j+1,k}^{n}+p_{i,j,k}^{n})\\& &-\theta
\frac{\Delta x\Delta z}{4a}(p_{i,j+1,k}^{n}-2p_{i,j,k}^{n}+p_{i,j-1,k}^{n}),
\end{eqnarray*}

\begin{eqnarray*}
F_{w _{i,j,k+\frac{1}{2}}^{n+1}}&=&\frac{w_{i,j,k+1}^{n+1}+w_{i,j,k}^{n+1}}{2}+\theta\frac{\Delta x\Delta y}{4a}(p_{i,j,k+2}^{n}-2p_{i,j,k+1}^{n}+p_{i,j,k}^{n})\\& &-
\theta\frac{\Delta x\Delta y}{4a}(p_{i,j,k+1}^{n}-2p_{i,j,k}^{n}+p_{i,j,k-1}^{n}),
\end{eqnarray*}
 $$\forall~i=0,\ldots,M, \forall~j=0,\ldots,N,
\forall~k=0,\ldots,L,$$ where: $\theta$
 is the relaxation
coefficient and
$$a=\frac{3 \Delta x \Delta y \Delta z}{2\Delta
t}+2\varepsilon\frac{\Delta x\Delta y}{\Delta
z}+2\varepsilon\frac{\Delta y\Delta z}{\Delta
x}+2\varepsilon\frac{\Delta x\Delta z}{\Delta y}.$$

\section{New finite volume discretization}
In this section we introduce a new finite volume schemes, that is we
approximate the solution of $(\ref{fluide1})$ by:
$$\boldsymbol{u}_{h}=\sum_{i,j=1}\boldsymbol{r}_{i,j,0}\boldsymbol{\hat{\overline{\varphi}}}^{0,\varepsilon}\chi_{i,j,0}
+\sum_{i,j=1}\boldsymbol{r_{i,j,L+1}}
\boldsymbol{\hat{\widetilde{\varphi}}}^{1,\varepsilon}\chi_{i,j,L+1}
+\sum_{i,j,1}\boldsymbol{u}_{i,j,k}\chi_{i,j,k},$$
where:
$$h=\Delta z,$$
$$\boldsymbol{r}_{i,j,0}=\frac{\boldsymbol{u}_{i,j,0}+\boldsymbol{u}_{i,j,1}}{2},$$
$$\boldsymbol{r}_{i,j,L+1}=\frac{\boldsymbol{u}_{i,j,L+1}+\boldsymbol{u}_{i,j,L}}{2},$$
$$\chi_{i,j,0}=\chi_{(x_{i-\frac{1}{2}},x_{i+\frac{1}{2}})\times(y_{j-\frac{1}{2}},y_{j+\frac{1}{2}})\times(0,h)},$$
$$\chi_{i,j,N+1}=\chi_{(x_{i-\frac{1}{2}},x_{i+\frac{1}{2}})\times(y_{j-\frac{1}{2}},y_{j+\frac{1}{2}})\times((L-1)h,Lh)},$$
$$\chi_{i,j,k}=\chi_{(x_{i-\frac{1}{2}},x_{i+\frac{1}{2}})\times(y_{j-\frac{1}{2}},y_{j+\frac{1}{2}})\times(z_{k-\frac{1}{2}},z_{k+\frac{1}{2}})},$$
and
\begin{eqnarray*}
{\hat{\overline{\varphi}}_{i}^{0,\varepsilon}}&=& -\int_{0}^{t}\frac{1}{\sqrt{4\pi(t-\tau)}} \frac{z}{2\sqrt{\varepsilon}(t-\tau)} e^{\frac{-z^2}{4\varepsilon (t-\tau)}}\\& &
\times\{2 \tau\cos(\alpha(\tau-t))-2\tau \sin(\alpha(\tau-t))\}d\tau,~\forall~i=1,2,
\end{eqnarray*}
\begin{eqnarray*}
{\hat{\widetilde{\varphi}}_{i}^{1,\varepsilon}}&=& -\int_{0}^{t}\frac{1}{\sqrt{4\pi(t-\tau)}} \frac{z}{2\sqrt{\varepsilon}(t-\tau)}
e^{\frac{-(1-z^2)}{4\varepsilon (t-\tau)}}\\& &
\times\{2 \tau\cos(\alpha(\tau-t))-2 \tau\sin(\alpha(\tau-t))\}d\tau,~\forall~i=1,2.
\end{eqnarray*}
\begin{eqnarray*}
{\hat{\overline{\varphi}}_{3}^{0,\varepsilon}}={\hat{\widetilde{\varphi}}_{3}^{1,\varepsilon}}=0.
\end{eqnarray*}

Multiplying $(\ref{fluide1})_{1}$ by $\chi_{i,j,k}$, integrating over $\Omega$, and replacing $\boldsymbol{u}^{\varepsilon}$ by $\boldsymbol{u}_{h}$ we find that the equations are the same as the classical finite volume scheme $(\ref{tem})$. Moreover the correctors verify $(\ref{fluide 4})_{1}$. Hence they do not contribute to these equations.
For the numerical simulations we do not use the modified boundary layer
$\boldsymbol{{\overline{\varphi} }}^{0,\varepsilon}$ and $\boldsymbol{{\widetilde{\varphi} }}^{1,\varepsilon}$ directly. Instead we consider another approximate form which reads as follows:
$$\boldsymbol{\widetilde{\overline{\varphi}}} ^{0,\varepsilon}(t,z)=(-exp(\frac{-z^2}{4\varepsilon t}),-exp(\frac{-z^2}{4\varepsilon t}),0).$$
Indeed, The approximation $\boldsymbol{\widetilde{\overline{\varphi}}} ^{0,\varepsilon}$ is much easier to be implemented numerically in coding than the theoretical corrector
$\boldsymbol{{\overline{\varphi}}^{0,\varepsilon}}$ obtained in section 2.\\
\newpage
Due to the nodes $\boldsymbol{r}_{i,j,0}$ and $\boldsymbol{r}_{i,j,L+1}$, the linear system associated with this scheme is not closed. However, by adding the correctors, we ensuring the closure of the linear system corresponding to the NFVM considered. Hence,
We multiply $(\ref{num})$ by the corrector
$\boldsymbol{\widetilde{\overline{\varphi}}}^{0,\varepsilon}$ and
integrate over $K_{i,j,1}$, we find:
\begin{eqnarray}\label{numer}
\int_{K_{ij1}}\!\!\!\!\dfrac{3\boldsymbol{u}^{n+1}-4\boldsymbol{u}^{n}+\boldsymbol{u}^{n-1}}{2\Delta
t}\boldsymbol{\widetilde{\overline{\varphi}}}
^{0,\varepsilon}-\varepsilon\int_{K_{ij1}}\!\!\Delta
\boldsymbol{u}^{n+1}\boldsymbol{\widetilde{\overline{\varphi}}}
^{0,\varepsilon}+\int_{K_{ij1}}\!\!\!\boldsymbol{\omega}\times
(2\boldsymbol{u}^{n}-\boldsymbol{u}^{n-1})\boldsymbol{\widetilde{\overline{\varphi}}}
^{0,\varepsilon}\\+2\int_{K_{ij1}}\nabla
p^{n}\widetilde{\overline{\varphi}}
^{0,\varepsilon}-\int_{K_{ij1}}\nabla
p^{n-1}\boldsymbol{\widetilde{\overline{\varphi}}}
^{0,\varepsilon}=\int_{K_{ij1}}\boldsymbol{f}^{n+1}\boldsymbol{\widetilde{\overline{\varphi}}}
^{0,\varepsilon}.\nonumber
\end{eqnarray}
In the following we calculate each term of $(\ref{numer})$, for the first term in the LHS (left-hand side) of $(\ref{numer})$
we find: $$ \int_{K_{ij1}}\dfrac{3\boldsymbol{u}^{n+1}-4\boldsymbol{u}^{n}+\boldsymbol{u}^{n-1}}{2\Delta
t}\boldsymbol{\widetilde{\overline{\varphi}}}
^{0,\varepsilon}dxdydz=\frac{3\boldsymbol{u}^{n+1}_{i,j,1}-4\boldsymbol{u}^{n}_{i,j,1}+\boldsymbol{u}^{n-1}_{i,j,1}}{2\Delta
t}\int_{K_{ij1}}\boldsymbol{\widetilde{\overline{\varphi}}}
^{0,\varepsilon}dxdydz.$$
For the second term in the LHS of $(\ref{numer})$, we obtain:
\begin{eqnarray}\label{done}
\int_{K_{ij1}}\Delta
\boldsymbol{u}^{n+1}\boldsymbol{\widetilde{\overline{\varphi}}}
^{0,\varepsilon}dxdydz&=&-\int_{K_{ij1}}\nabla
\boldsymbol{u}^{n+1}\nabla\boldsymbol{\widetilde{\overline{\varphi}}}
^{0,\varepsilon}dxdydz+\int_{\partial
K_{ij1}}\boldsymbol{\widetilde{\overline{\varphi}}} ^{0,\varepsilon}\frac{\partial\boldsymbol{u}^{n+1}}{\partial n}d\Gamma,\\
\nonumber&=&-\int_{K_{ij1}}\frac{\partial \boldsymbol{u}^{n+1}}{\partial
z}\frac{\partial\boldsymbol{\widetilde{\overline{\varphi}}}
^{0,\varepsilon}}{\partial z}dxdydz+\int_{\partial
K_{ij1}}\boldsymbol{\widetilde{\overline{\varphi}}}
^{0,\varepsilon}\frac{\partial \boldsymbol{u}^{n+1}}{\partial n}d\Gamma.\nonumber
\end{eqnarray}
Now, we calculate the first term in the RHS of $(\ref{done})$ we find:
\begin{eqnarray*}
\int_{K_{ij1}}\nabla \boldsymbol{u}^{n+1}\nabla\boldsymbol{\widetilde{\overline{\varphi}}^{0,\varepsilon}}
dxdydz&=&\int_{K_{ij1}}\frac{\partial\boldsymbol{u}^{n+1}}{\partial z}
\frac{\partial\boldsymbol{\widetilde{\overline{\varphi}}}^{0,\varepsilon}}{\partial z}dxdydz
\\&=&
\int_{x_{i-1/2}}^{x_{i+1/2}}\int_{y_{i-1/2}}^{y_{i+1/2}}\int_{0}^{h/2}\frac{\partial
\boldsymbol{u}^{n+1}}{\partial z}\frac{\partial\boldsymbol{\widetilde{\overline{\varphi}}}
^{0,\varepsilon}}{\partial z}dxdydz\\& &
+\int_{x_{i-1/2}}^{x_{i+1/2}}\int_{y_{i-1/2}}^{y_{i+1/2}}\int_{h/2}^{h}\frac{\partial
\boldsymbol{u}^{n+1}}{\partial z}\frac{\partial\boldsymbol{\widetilde{\overline{\varphi}}}
^{0,\varepsilon}}{\partial z} dxdydz\\&=&\frac{\boldsymbol{u}^{n+1}_{i,j,1}-\boldsymbol{r}^{n+1}_{i,j,1}}{h/2}\int_{x_{i-1/2}}^{x_{i+1/2}}\int_{y_{i-1/2}}^{y_{i+1/2}}
\int_{0}^{h/2}\frac{\partial\boldsymbol{\widetilde{\overline{\varphi}}}
^{0,\varepsilon}}{\partial z}dxdydz\\&
&+\frac{\boldsymbol{u}^{n+1}_{i,j,2}-\boldsymbol{u}^{n+1}_{i,j,1}}{h}\int_{x_{i-1/2}}^{x_{i+1/2}}\int_{y_{i-1/2}}^{y_{i+1/2}}\int_{h/2}^{h}
\frac{\partial\boldsymbol{\widetilde{\overline{\varphi}}}
^{0,\varepsilon}}{\partial z}dxdydz\\&=&
\frac{2}{h}(\boldsymbol{u}^{n+1}_{ij1}-\boldsymbol{r}^{n+1}_{ij0})\Delta x\Delta
y(\boldsymbol{\widetilde{\overline{\varphi}}}
^{0,\varepsilon}(\frac{h}{2})-\boldsymbol{\widetilde{\overline{\varphi}}}
^{0,\varepsilon}(0))\\&
&+\frac{\boldsymbol{u}^{n+1}_{i,j,2}-\boldsymbol{u}^{n+1}_{i,j,1}}{h} \Delta x\Delta
y(\boldsymbol{\widetilde{\overline{\varphi}}}
^{0,\varepsilon}(h)-\boldsymbol{\widetilde{\overline{\varphi}}}
^{0,\varepsilon}(h/2)).
\end{eqnarray*}
For the second term in the RHS of $(\ref{done})$ we obtain:
\begin{eqnarray*} \int_{\partial
K_{ij1}}\boldsymbol{\widetilde{\overline{\varphi}}}
^{0,\varepsilon}\frac{\partial \boldsymbol{u}^{n+1}}{\partial
n}d\Gamma&=&\int_{x_{i-1/2}}^{x_{i+1/2}}\int_{y_{i-1/2}}^{y_{i+1/2}}
\int_{z=0}\boldsymbol{\widetilde{\overline{\varphi}}}
^{0,\varepsilon}(-\frac{\partial \boldsymbol{u}}{\partial
z})d\Gamma\\& &+\int_{x_{i-1/2}}^{x_{i+1/2}}\int_{y_{i-1/2}}^{y_{i+1/2}}
\int_{z=h}\boldsymbol{\widetilde{\overline{\varphi}}}
^{0,\varepsilon}(\frac{\partial \boldsymbol{u}}{\partial
z})d\Gamma\\& &+\int_{x_{i-1/2}}^{x_{i+1/2}}\int_{0}^{h}
\int_{y=y_{j-1/2}}\boldsymbol{\widetilde{\overline{\varphi}}}
^{0,\varepsilon}(-\frac{\partial \boldsymbol{u}}{\partial
y})d\Gamma\\& &+\int_{x_{i-1/2}}^{x_{i+1/2}}\int_{0}^{h}
\int_{y=y_{j+1/2}}\boldsymbol{\widetilde{\overline{\varphi}}}
^{0,\varepsilon}(\frac{\partial \boldsymbol{u}}{\partial y})d\Gamma
\\& &+\int_{y_{i-1/2}}^{y_{i+1/2}}\int_{0}^{h}
\int_{x=x_{i-1/2}}\boldsymbol{\widetilde{\overline{\varphi}}}
^{0,\varepsilon}(-\frac{\partial \boldsymbol{u}}{\partial
x})d\Gamma\\& &+\int_{y_{i-1/2}}^{y_{i+1/2}}\int_{0}^{h}
\int_{x=x_{i+1/2}}\boldsymbol{\widetilde{\overline{\varphi}}}
^{0,\varepsilon}(\frac{\partial \boldsymbol{u}}{\partial x})d\Gamma.
\end{eqnarray*}
Now, the third term in the LHS of $(\ref{numer})$, can be rewritten as bellow:
\begin{eqnarray*}
\int_{K_{ij1}}\boldsymbol{\omega}\times
(2 \boldsymbol{u}^{n}-\boldsymbol{u}^{n-1})\boldsymbol{\widetilde{\overline{\varphi}}}
^{0,\varepsilon}dxdydz=\boldsymbol{\omega}\times (2 \boldsymbol{u}^{n}_{i,j,1}-\boldsymbol{u}^{n-1}_{i,j,1})\Delta x \Delta y \int_{0}^{h}\boldsymbol{\widetilde{\overline{\varphi}}}^{0,\varepsilon}dxdydz.
\end{eqnarray*}
We calculate the first component of the fourth term in the LHS of $(\ref{numer})$ and, we find:
\begin{eqnarray*} \int_{K_{ij1}}\partial_{x}
p^{n+1}\widetilde{\overline{\varphi}}_{1}
^{0,\varepsilon}dx dy dz=\frac{p^{n+1}_{i+1,j,1}-p^{n+1}_{i-1,j,1}}{2\Delta
x}\Delta x \Delta y \int_{0}^{h}\widetilde{\overline{\varphi}}_{1}
^{0,\varepsilon} dz.
\end{eqnarray*}
For the second component of the fourth term in the LHS of $(\ref{numer})$ we have:
\begin{eqnarray*} \int_{K_{ij1}}\partial_{y}
p^{n+1}\widetilde{\overline{\varphi}}_{2}
^{0,\varepsilon}dx dy dz=\frac{p^{n+1}_{i,j+1,1}-p^{n+1}_{i,j-1,1}}{2\Delta
y}\Delta x \Delta y \int_{0}^{h}\widetilde{\overline{\varphi}}_{2}
^{0,\varepsilon} dz.
\end{eqnarray*}
Concerning the first term on the RHS of $(\ref{numer})$, we have:
\begin{eqnarray*}
\int_{K_{ij1}}\boldsymbol{f}^{n+1}\boldsymbol{\widetilde{\overline{\varphi}}}
^{0,\varepsilon}dxdydz=\Delta x\Delta y \boldsymbol{f}^{n+1}_{i,j,1}\int_{0}^{h}\boldsymbol{\widetilde{\overline{\varphi}}}
^{0,\varepsilon}dz.
\end{eqnarray*}
Hence, we infer that
\begin{eqnarray*}
& &\frac{1}{2\Delta
t}(3\boldsymbol{u}^{n+1}_{i,j,1}-4\boldsymbol{u}^{n}_{i,j,1}+\boldsymbol{u}^{n-1}_{i,j,1})\int_{0}^{h}\boldsymbol{\widetilde{\overline{\varphi}}}
^{0,\varepsilon}dz
-\varepsilon [\frac{1}{h}(-3\boldsymbol{\widetilde{\overline{\varphi}}}
^{0,\varepsilon}(h/2)\boldsymbol{u}^{n+1}_{i,j,1}+\\& &+2\boldsymbol{r}^{n+1}_{i,j,0}\boldsymbol{\widetilde{\overline{\varphi}}}
^{0,\varepsilon}(h/2)
+\boldsymbol{u}^{n+1}_{i,j,2}\boldsymbol{\widetilde{\overline{\varphi}}}
^{0,\varepsilon}(h/2))
-\big(\frac{1}{(\Delta x)^2}(\boldsymbol{u}^{n+1}_{i-1,j,1}-2\boldsymbol{u}^{n+1}_{i,j,1}+\boldsymbol{u}^{n+1}_{i+1,j,1})+\\& &
+\frac{1}{(\Delta y) ^2}(\boldsymbol{u}^{n+1}_{i,j-1,1}
-2 \boldsymbol{u}^{n+1}_{i,j,1}+\boldsymbol{u}^{n+1}_{i,j+1,1})\big)\int_{0}^{h}\boldsymbol{\widetilde{\overline{\varphi}}}
^{0,\varepsilon}dz]+\boldsymbol{\omega}\times(2 \boldsymbol{u}_{i,j,1}^{n}-\boldsymbol{u}_{i,j,1}^{n-1}) \int_{0}^{h}\boldsymbol{\widetilde{\overline{\varphi}}}^{0,\varepsilon}dz+\\& &
+2\left(
    \begin{array}{c}
(\frac{p^{n}_{i+1,j,1}-p^{n}_{i-1,j,1}}{2\Delta x})\int_{0}^{h}\boldsymbol{\widetilde{\overline{\varphi}}} ^{0,\varepsilon}dz \\\\
(\frac{p^{n}_{i,j+1,1}-p^{n}_{i,j-1,1}}{2\Delta y})\int_{0}^{h}\boldsymbol{\widetilde{\overline{\varphi}}} ^{0,\varepsilon}dz   \\\\
      0\\
    \end{array}
  \right)
-
\left(
    \begin{array}{c}
(\frac{p^{n-1}_{i+1,j,1}-p^{n-1}_{i-1,j,1}}{2\Delta x})\int_{0}^{h}\boldsymbol{\widetilde{\overline{\varphi}}} ^{0,\varepsilon}dz \\\\
(\frac{p^{n-1}_{i,j+1,1}-p^{n-1}_{i,j-1,1}}{2\Delta y})\int_{0}^{h}\boldsymbol{\widetilde{\overline{\varphi}}} ^{0,\varepsilon}dz   \\\\
      0\\
    \end{array}
  \right)
\\& &=\boldsymbol{f}^{n+1}_{i,j,1}\int_{0}^{h}\boldsymbol{\widetilde{\overline{\varphi}}}
^{0,\varepsilon}dz.
\end{eqnarray*}

\section{Numerical results} In this section we will
compute the error approximation using the classical finite volume
method and the new finite volume method, so the pressure and the
source term are chosen such that:

\begin{equation*}
p(x,y,z,t)=\cos(2\pi x ) \cos(2\pi y)cos(\pi z) t,
\end{equation*}

\begin{eqnarray*}
 u^{\varepsilon}(x,y,z,t)=t \sin(2\pi y)
 (1-e^{\frac{-z}{\sqrt{\varepsilon}}}\cos(\frac{z}{\sqrt{\varepsilon}})) (1-e^{\frac{-(1-z)}{\sqrt{\varepsilon}}}\cos(\frac{(1-z)}{\sqrt{\varepsilon}})),
\end{eqnarray*}
\begin{eqnarray*}
 v^{\varepsilon}(x,y,z,t)=t \sin(2\pi x)
 (1-e^{\frac{-z}{\sqrt{\varepsilon}}}\cos(\frac{z}{\sqrt{\varepsilon}})) (1-e^{\frac{-(1-z)}{\sqrt{\varepsilon}}}\cos(\frac{(1-z)}{\sqrt{\varepsilon}})),
\end{eqnarray*}
and
\begin{eqnarray*}
 w^{\varepsilon}(x,y,z,t)=0.
\end{eqnarray*}
Note that the test solution given above satisfies the equations $(\ref{fluide1})_{1,2}$ including the boundary and initial conditions $(\ref{fluide1})_{3,5}$  with $\boldsymbol{u}_{0}=0$.\\
Thus, the function source is chosen using this test solution.

\begin{figure}
\begin{center}
\begin{tabular}{|c|c|c|c|c|}   \hline
   N=M=L & t & $\varepsilon$ & CFVM & NFVM \\
   \hline
   10 & 1 & $10^{-2}$ & 0.03206& 0.12836 \\

   20& 1 & $10^{-2}$ & 0.00634 & 0.03893\\

   30& 1 &$10^{-2}$& 0.00269 &0.02553\\
   \hline
   10&1&$10^{-3}$&0.092294&0.22647\\

   20&1&$10^{-3}$& 0.033726&0.15753\\

   30&1&$10^{-3}$& 0.01331& 0.08020\\
   \hline
   10& 1 &$10^{-5}$& 1.61660e+03 &0.04487\\

   20& 1 &$10^{-5}$& 0.08741&0.010303\\

   30& 1 &$10^{-5}$& 0.11722 &0.00460\\
   \hline
   10&1&$10^{-6}$&1.10612e+10&0.044901\\
   20&1&$10^{-6}$&4.42881e+06&0.01032\\
   30&1&$10^{-6}$&1.12960e+03&0.00442\\
   \hline
   10&1&$10^{-7}$&5,26218e+62&0.04490\\
   20&1&$10^{-7}$&1.16428e+29&0.01032\\
   30&1&$10^{-7}$&6.72495e+17&0.00443\\
  \hline
\end{tabular}
 \caption{The $L^{2}$ norm of the velocity error with classical
finite volume (CFVM) and new finite volume method (NFVM) for different values of $\varepsilon$ at $t=1$}.
 \end{center}
\end{figure}

\begin{figure}
\begin{center}
\begin{tabular}{|c|c|c|c|c|}
  \hline
  N=M=L & t & $\varepsilon$ & CFVM & NFVM \\
  \hline
  10& 1 & $10^{-2}$ & 0.02493 &  0.03178\\

  20 & 1 & $10^{-2}$ & 0.00511& 0.00920 \\

  30 & 1 & $10^{-2}$ & 0.00224& 0.00533 \\
  \hline
  10&1&$10^{-3}$& 0.02684 &0.02771\\

  20&1&$10^{-3}$&0.00553& 0.00907\\
  30&1&$10^{-3}$&0.002381&0.00590\\
  \hline
  10 & 1 & $10^{-5}$ & 1.48996e+02 & 0.02602 \\

  20 & 1 & $10^{-5}$ & 0.00774& 0.00539 \\

  30 & 1 & $10^{-5}$ & 0.00655& 0.00238 \\
  \hline
  10&1&$10^{-6}$&1.01953e+16&0.026016\\
  20&1&$10^{-6}$&2.83861e+05&0.00539\\
  30&1&$10^{-6}$&58.98117&0.00238\\
  \hline
  10&1&$10^{-7}$&4.85027e+61&0.02601\\
  20&1&$10^{-7}$&7.46273e+27&0.005394\\
  30&1&$10^{-7}$&3.51186e+16&0.00238\\
  \hline
\end{tabular}
\caption{The $L^{2}$ norm of the pressure error with classical
finite volume (CFVM) and new
finite volume method (NFVM) for different values of $\varepsilon$ at $t=1$}.
\end{center}
\end{figure}

\newpage

\newpage
%
%

\section{Conclusion and fracture works}
In this paper we have compared two different finite volume methods (CFVM) and (NFVM) when the viscosity is considered small and more precisely
in the rang $10^{-3}-10^{-7}$ . We derived an approximate solution of the-time dependent rotating fluid in $3D$ channel using the splitting methods for the time discretization and colocated space discretization. One of the novelties of this article is that we propose a new numerical approach to treat the pressure and the divergence free condition introducing correctors to solve the boundary layers. We also show that the (NFVM) is more performing than the (CFVM) when the viscosity is small. We showed that our NFVM still perform for very large Reynolds number. To the best of our knowledge, this is the first work which gives a new finite volume scheme taking into account boundary layer variations for the linearized Navier-Stokes equations. Note that the method developed here may apply to many other problems and domains. This will be the subject of subsequent work.

\section{Appendix.}
In this paragraph, we give a sketch of the proof of the existence and regularity of the
solution of the limit problem $(\ref{eq lim 00})$. For complete study of the existence of system $(\ref{eq lim 00})$ we refer the reader to , see also     and      .We first want to apply the Hille-Phillips-Yosida Theorem to prove the existence and uniqueness of the
solution of $(\ref {eq lim 00})$. Thus we start by introducing
the adequate function spaces:
\begin{eqnarray*}
H &=& \{\boldsymbol{v}\in(L^{2}(\Omega))^{3};\mbox{div}\boldsymbol{v}=0, v_{3}(z=0)=v_{3}(z=h)=0, \\
   &&\mbox{and $\boldsymbol{v}$ is $2\pi$ periodic in the $x$ and $y$ directions}\}.
\end{eqnarray*}
\begin{equation}
D(A)=\{\boldsymbol{v}\in H;\exists p\in \mathbf{D'}(\Omega),~\text{such}~
\boldsymbol{\omega}\times \boldsymbol{v}+\nabla p\in H\},\nonumber
\end{equation}
with the norm:
\begin{equation}
\|\boldsymbol{v}\|_{D(A)}=(\|\boldsymbol{v}\|_{H}^{2}+\|\boldsymbol{\omega}\times \boldsymbol{v}+\nabla
p\|_{H}^{2})^{1/2}.
\end{equation}
Then for $\boldsymbol{v}\in D(A)$ we set $A\boldsymbol{v}=\boldsymbol{\omega}\times \boldsymbol{v}+\nabla
p$, thus we define an unbounded linear operator $A$ which maps $D(A)\subset
H $onto  $H$.
\begin{theorem}[Hille-Yosida Theorem]
Let $H$ be a Hilbert space and let $B:D(B)\longrightarrow H$
 a linear unbounded operator, with domain
$D(B)\subset H$ such that $D(B)$ is dense in $H$ and $(-B)$ is
m-dissipative. Then $(-B)$ is the infinitesimal generator of a
contraction semigroup $\{S(t)\}_{t>0}$ in $H$, and the solution of
the following system:
\begin{equation}\label{eq lim 0}
\left\lbrace
\begin{array}{lcl}
\dfrac{d\boldsymbol{v}}{dt}+B\boldsymbol{v}=\boldsymbol{f},\\
\boldsymbol{v}|_{t=0}=\boldsymbol{v}_{0},
\end{array}\right.
\end{equation}
satisfies the following properties:\\
$(H_{0})$ If $\boldsymbol{v}_{0}$ and $f\in L^{1}(0,T;H)$, then $\boldsymbol{v}\in C([0,T];H), \forall~T>0$.\\
$(H_{1})$ If $\boldsymbol{v}_{0}\in D(B)$ and $f'\in L^{1}(0,T;H)$ then $ \boldsymbol{v}\in
C^{1}([0,T];H)\cap C^{0}([0,T];D(B))$ and $\dfrac{d\boldsymbol{v}}{dt}\in
L^{\infty}([0,T];H), \forall~T>0.$
\end{theorem}

\begin{remark}
A linear operator is dissipative in $H$ if and only if:  $\forall
\boldsymbol{u}\in D(A), \forall\lambda>0, \|\boldsymbol{u}-\lambda A\boldsymbol{u}\|\geq\|\boldsymbol{u}\|.$
\end{remark}

\begin{remark}
A linear operator $A$ is m-dissipative if: A is dissipative and
$\forall f\in X, \forall\lambda>0, \exists \boldsymbol{u}\in D(A), \boldsymbol{u}-\lambda
A\boldsymbol{u}=\boldsymbol{f}.$
\end{remark}

\begin{preuve}\normalfont Now we want to show that the operator (-A) is
m-dissipative, hence we will prove that the following system:
\begin{equation}\label{yosida}
\left\lbrace
\begin{array}{lcl}
\lambda\boldsymbol{\omega}\times \boldsymbol{u}+\lambda \nabla p+\boldsymbol{u}&=&\boldsymbol{f},\\
$div$~\boldsymbol{u}=0,\\
u_{3}=0,~\mbox{en}\hspace{0,2cm} z=0,1,
\end{array}\right.
\end{equation}
has a unique solution in $D(A)$ for all   $\boldsymbol{f}\in H$ and $\forall
\lambda>0$, and the solution satisfies the estimate:
\begin{equation}
\|\boldsymbol{u}\|_{H}\leq \|\boldsymbol{f}\|_{H},\hspace{0,5 cm} \forall \boldsymbol{f}\in H.
\end{equation}
We multiply $(\ref{yosida})$ by $\boldsymbol{v}\in H$, integrate over
$\Omega$ and we find:
$$\lambda\int_{\Omega}(\boldsymbol{\omega}\times \boldsymbol{u}).\boldsymbol{v} d\Omega +\lambda\int_{\Omega}\nabla p.\boldsymbol{v} d\Omega+\int_{\Omega}\boldsymbol{u}.\boldsymbol{v} d\Omega=\int_{\Omega}\boldsymbol{f}\boldsymbol{v} d\Omega.$$
We have:
\begin{eqnarray*}
\int_{\Omega}\nabla p \boldsymbol{v}&=&-\int_{\Omega}p\textnormal{div}\boldsymbol{v}+\int_{\partial\Omega}p \boldsymbol{v}.n d(\Gamma)=0.\\
\end{eqnarray*}
We set
$$a(\boldsymbol{u},\boldsymbol{v})=\lambda\int_{\Omega}(\boldsymbol{\omega}\times \boldsymbol{u}).\boldsymbol{v}+\int_{\Omega}\boldsymbol{u}.\boldsymbol{v},$$
and
$$F(\boldsymbol{v})=\int_{\Omega}\boldsymbol{f}\boldsymbol{v}.$$
Here $a(.,.)$ is a continuous and coercive bilinear form in
$H\times H$. In fact we have:
\begin{eqnarray*}
|a(\boldsymbol{u},\boldsymbol{v})|&\leq& \lambda|\boldsymbol{u}|_{H}|\boldsymbol{v}|_{H}+|\boldsymbol{u}|_{H}|\boldsymbol{v}|_{H},\\
&\leq&k(\lambda)|u|_{H}|v|_{H},
\end{eqnarray*}
and
$$|a(\boldsymbol{u},\boldsymbol{u})|=|\boldsymbol{u}|_{H}^{2}.$$
Also $F(\boldsymbol{v})$ is a continuous linear form:
$$\int_{\Omega}\boldsymbol{f}\boldsymbol{v}\leq |\boldsymbol{f}||\boldsymbol{v}|.$$
Hence according to the Lax-Milligram theorem, there exists a unique
$\boldsymbol{u}\in H$ such that:
$$\lambda~A{u}+\boldsymbol{u}=\boldsymbol{f},$$
that is,
\begin{eqnarray*}
\lambda \boldsymbol{\omega}\times \boldsymbol{u}+\lambda\nabla p+\boldsymbol{u}&=&\boldsymbol{f}.
\end{eqnarray*}
Multiplying the above equation by $\boldsymbol{u}$ and integrating over
$\Omega$, we find:
\begin{eqnarray*}
\lambda\int_{\Omega}(\boldsymbol{\omega}\times \boldsymbol{u})
\boldsymbol{u}+\lambda\int_{\Omega}\nabla p \boldsymbol{u}+\int_{\Omega}\boldsymbol{u} \boldsymbol{u}=\int_{\Omega}\boldsymbol{f}
\boldsymbol{u},
\end{eqnarray*}
then the solution $\boldsymbol{u}$ satisfies the estimate:
\begin{eqnarray*}
\|\boldsymbol{u}\|_{H}\leq \|\boldsymbol{f}\|_{H}.
\end{eqnarray*}
Also we have:
\begin{eqnarray*}
\|\boldsymbol{u}\|_{D(A)}&=&(\|\boldsymbol{u}\|_{H}^{2}+\|\omega\times \boldsymbol{u}+\nabla p\|_{H}^{2})^{1/2},\\
&\leq& \|\boldsymbol{u}\|_{H}+\|\boldsymbol{\omega}\times \boldsymbol{u}+\nabla p\|_{H},\\
&\leq& k(\lambda) \|\boldsymbol{f}\|_{H}.
\end{eqnarray*}
Hence (-A) is m-dissipative operator, Moreover we have $\boldsymbol{u}_{0}\in H$,
then according to the Hille-Yosida theorem the system $(\ref{eq
lim 00})$ has a unique solution $\boldsymbol{u}\in C([0,\infty[,H).$ Furthermore, we
have:
\begin{eqnarray*}
\|\nabla p\|_{H^{-1}}&\leq&\frac{1}{\lambda}\|\boldsymbol{f}\|_{H^{-1}}+\frac{1}{\lambda}\|\boldsymbol{u}\|_{H^{-1}}+\|\boldsymbol{\omega}\times \boldsymbol{u}\|_{H^{-1}}\\
&\leq&k(\lambda) \|\boldsymbol{f}\|_{H}.
\end{eqnarray*}
Then, we obtain;
\begin{equation*}
\|p\|_{L^{2}(\Omega)}\leq k(\lambda) \|\boldsymbol{f}\|_{H}.
\end{equation*}
\end{preuve}
\section{Acknowledgements.}
We are very grateful to Sylvain Faure for his collaboration to programming the code Matlab in section $4$.

\end{document}